\documentclass[11pt,a4paper]{article}
\usepackage{tikz,pgfplots}
\usepackage[]{times}
\usepackage{fullpage}
\usepackage{graphicx}
\usepackage{subfigure}
\usepackage{amsmath}
\usepackage{fancyhdr}
\date{}

\newcommand{\bw}{{\bf w}}
\newtheorem{theorem}{Theorem}[section]

\newtheorem{e-proposition}[theorem]{Proposition}

\newtheorem{e-definition}[theorem]{Definition\rm}
\newtheorem{remark}{\it Remark\/}

\title{Hermite finite elements for convection-diffusion equations}
  \author{
     F. A. Radu$^{1}$\\[1mm]
  {\small $^{1}$ Department of Mathematics, University of Bergen, Bergen, Norway}\\
  {\small e-mail: {\it florin.radu@math.uib.no}}\\[1mm]
    V. Ruas$^{2,3}$\\[1mm]
  {\small $^{2}$ Sorbonne Universit\'es, UPMC Univ Paris 06 \& CNRS, UMR 7190, IJRDA, 75005 Paris, France}\\
  {\small e-mail: {\it vitoriano.ruas@upmc.fr}}\\[1mm]
  {\small $^{3}$ CNPq scholar at Graduate school of Metrology for Quality and Innovation, PUC-Rio, Rio de Janeiro, Brazil}\\
  {\small e-mail: {\it vitoriano.ruas@pq.cnpq.fr}}\\[1mm]
   P. R. Trales$^{4}$\\[1mm]
  {\small $^{4}$ Department of Analysis, IME-Universidade Federal Fluminense, Niter\'oi, Brazil}\\
  {\small e-mail: {\it paulotrales@im.uff.br}}}
  
\begin{document}
\maketitle
\thispagestyle{fancy}

\begin{abstract}
This work addresses techniques to solve convection-diffusion problems based on Hermite interpolation. 
More specifically we extend to the case of these equations a Hermite finite element method  
providing flux continuity across inter-element boundaries, shown to be a well-adapted tool for simulating pure 
diffusion phenomena \cite{JCP}. We consider two methods that can be viewed as non trivial improved versions of the lowest order 
Raviart-Thomas mixed method \cite{RaviartThomas}, corresponding to its extensions to convection-diffusion problems proposed by 
Douglas and Roberts \cite{Douglas}. A detailed convergence study is carried out for one of the methods, and numerical 
results illustrate the performance of both of them, as compared to each other and to the corresponding mixed methods.
\end{abstract}


\maketitle

\section{Introduction}

\hspace{4mm} Historically Hermite finite elements have mostly been used to solve fourth order partial differential equations, 
because minimum continuity of solution derivatives across inter-element boundaries is required in this case. However the 
construction of such elements can be rather laborious, as shown in \cite{Ciarlet}. It is noticeable in this respect that the 
recent technique of the virtual element led to feasible constructions of $C^k$ functions for $k \geq 1$ on meshes consisting of 
straight elements of arbitrary shape \cite{virtual}, though by means of polynomials of rather high degree.\\  
\indent On the other hand Hermite interpolation has been showing to be a good alternative to solve several kinds of field problems 
modeled by second order boundary value problems in many respects. An outstanding demonstration of such an assertion is provided by 
the isogeometric analysis (IGA) introduced about ten years ago (see e.g. \cite{Hughes}). In this case advantage is taken from data 
satisfying high continuity requirements supplied by CAD, for a subsequent finite element analysis. However IGA in connection with 
triangular or tetrahedral meshes is incipient, in spite of the undeniable geometric flexibility of this type of partitions.
This is a good reason to study Hermite finite elements methods defined upon triangles or tetrahedra to solve second order partial 
differential equations, which are low order and easy to implement at a time. That is what we do in this work, by focusing more 
particularly the representation of fluxes for the simulation of phenomena or processes of the convection-diffusion type.  \\
\indent In practice quantities directly depending on partial derivatives of the variable in terms of which an equation is expressed, 
i.e., the primal variable, are often more important than this unknown itself. Among them one might quote the flux in a porous 
medium flow or in heat flow. As far as methods allowing to enforce the continuity of normal derivatives or normal fluxes across the 
boundaries of triangular or tetrahedral cells are concerned, both mixed finite elements and finite volumes have been playing a 
prominent role since long. However Hermite interpolation can also be a tool well adapted for this purpose, as shown in \cite{JCP} 
for pure diffusion equations in highly heterogeneous media. In this work and in \cite{CAM2013} two finite element methods of the 
Hermite type based on a quadratic interpolation were studied. Both have, either identical or better convergence properties than 
some of the classical mixed methods, such as the $RT_0$ method, i.e. the lowest order Raviart-Thomas' \cite{RaviartThomas}, 
according to the norm under consideration, though at comparable implementation cost. \\
\indent As far as the convection-diffusion equations are concerned, a rather great amount of numerical solution techniques are 
available today. Nevertheless the fact that these equations lie on the basis of the mathematical modeling  
of countless physical phenomena, keeps encouraging specialists in the search for 
efficient methodology to solve this class of problems. This is particularly true of 
convection dominated processes, in which the correct capture of sharp boundary layers often reveal demerits of widespread 
computational techniques, even when the problem to solve is linear.\\ 
\indent The main purpose of this work is to carry out a complete mathematical study of the Hermite method to solve the convection-
diffusion equations introduced in \cite{ICNAAM2013}. More specifically such a method is an extension to the convection-diffusion 
equations of the Hermite finite element studied in \cite{JCP} that can be regarded as a variant of the $RT_0$ mixed element 
\cite{CAM2013}. The method is described in Section 2, where we recall that it is uniformly stable with respect to a suitable working norm. In Section 3 we apply these results which immediately lead to a priori error estimates in the same norm. Naturally enough they are of the first order. However similarly to \cite{CAM2013} we also prove that the error in the $L^2$-norm is in 
terms of the square of the mesh size, in contrast to the first order ones that hold for the mixed extension of the $RT_0$ method to 
convection-diffusion-reaction equations in non divergence form proposed in \cite{Douglas}. In Section 4 we further consider a 
variant of the method under study that can be viewed as the Hermite analog of the Douglas \& Roberts mixed method \cite{Douglas}, 
applying to the C-D equations in divergence form. This variant is slightly different from the one introduced in \cite{HEFAT}. 
The convergence properties and the accuracy of both new methods are checked and compared by means of numerical experiments reported in 
Section 5. In Section 6 we conclude with some comments on the whole work. \\
\indent Referring to \cite{Adams}, in the sequel we employ the following notations: $S$ being a bounded open set of $\Re^N$, 
we denote the standard norm of Sobolev spaces $H^m(S)$ (resp. $W^{m,p}(S)$ for $p \geq 1$, $p \neq 2$),   
for any non negative integer $m$ by $\parallel \cdot \parallel_{m,S}$ (resp. $\parallel \cdot \parallel_{m,p,S}$), 
including $L^2(S)=H^0(S)$. The standard semi-norm of $H^m(\Omega)$ (resp. $W^{m,p}(S)$ for $p \geq 1$, $p \neq 2$) is denoted by 
$| \cdot |_{m,S}$ (resp. $| \cdot |_{m,p,S}$).
Finally $\forall S \subset \Omega$ and $\forall f,g \in L^2(S)$, $(f,g)_S:=\int_S f g \; dS$ if $S \neq \Omega$ and 
$(f,g):=\int_{\Omega} f g \; dx$.  \\
\indent Let $\Omega$ be a bounded domain of $\Re^N$, $N=2,3$, with boundary $\Gamma$, $f \in L^2(\Omega)$ be given a function 
$f \in L^2(\Omega)$, ${\cal K}$ be
a tensor assumed to be constant, symmetric and positive-definite and  let ${\bf w}\in [C^{0}(\bar{\Omega})]^N$ denote a velocity 
field. In this work we study as a model the following equation, assumed to have a unique solution:\\

Find $u$ such that $u=0 $ on $\Gamma$ and\\
\begin{equation}
\label{model}
 - \nabla \cdot {\cal K} \nabla u + {\bf w} \cdot \nabla u = f \mbox{ in } \Omega.
\end{equation}


\section{A Hermite solution method}

\hspace{4mm} Henceforth we assume that $\Omega$ is a polygon if $N=2$ or a polyhedron if $N=3$, and that 
we are given a finite element partition ${\cal T}_h$ of $\Omega$, consisting of triangles or tetrahedra according to 
the value of $N$, and belonging to a regular family of partitions (cf. \cite{Ciarlet}). 
$h$ denotes the maximum diameter of the elements of ${\cal T}_h$. \\
\indent In the following we define two finite element spaces $U_h$ and $V_h$ associated with ${\mathcal T}_h$.  
 Let ${\bf w}_h$ be the constant field in each element of $T \in {\mathcal T}_h$ whose value in $T$ is ${\bf w}({\bf x}_T)$, where 
${\bf x}_T$ is the position vector of the centroid of $T$, and ${\bf w}^1_h$ be the standard continuous piecewise linear 
interpolate of ${\bf w}$ at the vertices of ${\mathcal T}_h$. 
We further introduce the operators $\Pi_T : L^2(T) \longrightarrow L^2(T)$ given by 
$\Pi_T[v]:= \int_T v dx / meas(T)$ for $T \in {\mathcal T}_h$, and $\Pi_h : L^2(\Omega) \longrightarrow L^2(\Omega)$ by 
$\Pi_h[v]_{|T} = \Pi_T[v_{|T}] \; \forall T \in {\mathcal T}_h$. Now throughout this work we will work with the following \\

\noindent \underline{Local algebraic structure of the Hermite finite element spaces:} \\
Every function $v \in V_h$ (resp. $\in U_h$) is such that in each element $T \in {\cal T}_h$ it is expressed by 
\begin{equation}
\label{v}
v_{|T}={\bf x}^{\mbox{t}}\{{\cal K}^{-1}[a {\bf x} /2 + {\bf b}]\} + d,
\end{equation} 
\noindent where ${\bf x}$ represents the space variable, ${\bf b}$ is a constant vector of $\Re^N$ and $a$ and $d$ are two real coefficients.\\
In every $N$-simplex $T$ we associate with a quadratic function $v$ of the form (\ref{v}):\\

\noindent \underline{Sets ${\mathcal D}_T$ and ${\mathcal E}_T$ of local degrees of freedom for the Hermite finite element spaces:} \\
$F$ being an edge if $N=2$ or a face if $N=3$ belonging to the boundary $\partial T$ of an $N$-simplex $T$, 
and ${\bf n}_F$ being the unit normal vector on $F$ oriented in a given manner for each $F \subset \partial T$, we set:
\begin{equation}
\label{U}
\left\{ 
\begin{array}{l}
{\mathcal D}_T:= \displaystyle \{\cup_{F} \; {\mathcal U}_F\} \cup {\mathcal U}_T \mbox{ where }\\
 {\mathcal U}_F(v)=\int_{F}{\cal K} \nabla v \cdot {\bf n}_F ds /meas(F); \\ 
{\mathcal U}_{T}(v) = \int_T v \;dx/meas(T).
\end{array} 
\right.
\end{equation}
\begin{equation}
\label{V} 
\left\{ 
\begin{array}{l}
{\mathcal E}_T:= \displaystyle \{ \cup_{F} \; {\mathcal V}_F \} \cup {\mathcal V}_T \mbox{ where }\\
 {\mathcal V}_F(v)=\int_{F} ({\cal K} \nabla v + {\bf w}_h \Pi_T[v] ) \cdot 
{\bf n}_F ds/meas(F); \\  
{\mathcal V}_{T}(v) = \int_T v \;dx/meas(T).
\end{array} 
\right.
\end{equation}
The canonical basis functions associated with these sets of degrees of freedom 
are as follows. First we note that $\forall v \in V_h$ or $\in U_h$, 
$\nabla v_{|T}$ for $T \in {\mathcal T}_h$ is expressed by 
${\cal K}^{-1} [a_T {\bf x} + {\bf b}_T]$ for certain $a_T \in \Re$ and ${\bf b}_T \in \Re^N$. 
Then the flux variable ${\mathcal K} \nabla v_{|T}$ is of the form $a_T {\bf x} + {\bf b}_T$, and from a well-known property of the lowest order Raviart-Thomas mixed element $a_T$ and ${\bf b}_T$ can be 
uniquely determined for prescribed ${\mathcal V}_F(v)$ (resp ${\mathcal U}_F(v)$), $\forall F \subset \partial T$. 
Indeed by construction the flux variable for the Hermite element is locally defined by functions of the same form as for the lowest order Raviart-Thomas element. Once $a_T$ and ${\bf b}_T$ are known, we 
determine the value of the additive constant $d_T$ to complete the expression of $v_{|T}$, by enforcing the condition $\Pi_T[v]=0$. 
As for the basis function corresponding to the degree of freedom ${\mathcal V}_{T}$ (resp. ${\mathcal U}_{T}$), 
the values of $a_T$ and ${\bf b}_T$ are obtained in a trivial manner as specified in \cite{ICNAAM2013}. Then the value of $d_T$ is adjusted in such a way that the mean value of the corresponding quadratic function is one. 
This should be enough to determine the $N+2$ basis functions associated with a given $N$-simplex $T$, corresponding to the sets 
${\mathcal U}_T$ and ${\mathcal V}_T$ of degrees of freedom, for spaces $U_h$ and $V_h$ respectively, 
since the $RT_0$ method is well-known (cf. \cite{Thomas}). However for the sake of clarity we exhibit them below.  \\
\noindent $T$ being an element of ${\mathcal T}_h$ let ${\bf x}_i^T$ be the position vector of the $i$-th vertex $S^T_i$ of $T$, 
$F^T_i$ be the face of $T$ opposite to $S^T_i$ and $h_i^T$ be the length of the corresponding height of $T$, for $i=1,\ldots,N+1$. 
We have:\\

\noindent \underline{Local basis functions $\varphi^T_i$ for space $U_h$:} \\
The local basis function $\varphi_i^T$ associated with the degree of freedom ${\mathcal U}_{F_i^T}$ and the basis function 
$\varphi_{N+2}^T$ associated with the degree of freedom ${\mathcal U}_T$ are given by:
\begin{equation} 
\label{basisU}
\left\{
\begin{array}{l}
\varphi_i^T = {\bf x}^{\mbox{t}}\{{\cal K}^{-1}[a_i^T {\bf x}/2 + {\bf b}_i^T]\} + d_i^T \mbox{ for } i=1,\ldots,N+2, \mbox{ where}\\
\left.
\begin{array}{l}
a_i^T=[h^T_i]^{-1}, \\ 
{\bf b}_i^T = -{\bf x}^T_i a_i^T, \\
d_i^T=-\int_T {\bf x}^{\mbox{t}}\{{\cal K}^{-1}[a_i^T {\bf x} /2 + {\bf b}_i^T]\} \; dx /meas(T), 
\end{array}
\right\}  
\mbox{for } i=1,\ldots,N+1;\\ 
\begin{array}{l}
a_{N+2}^T=0, \\ 
{\bf b}_{N+2}^T = (0,\ldots,0)^{\mbox{t}}, \\
d_{N+2}^T=1.
\end{array}
\end{array}
\right.
\end{equation}
\noindent \underline{Local basis functions $\psi^T_i$ for space $V_h$:} \\
Akin to the case of $U_h$, the $\psi^T_i$'s are functions of the form (\ref{v}) for $i=1,\ldots,N+2$. 
Since by definition the mean values of all the first $N+1$ $\varphi_i^T$'s vanish, the local basis functions  
$\psi_i^T$ of $V_h$ associated with the degree of freedom ${\mathcal V_{F_i^T}}$ for $i=1,\ldots,N+1$, together with its local basis function $\psi_{N+2}^T$ associated with the degree of freedom ${\mathcal V}_T$ are given by:  
\begin{equation} 
\label{basisV}
\left\{ 
\begin{array}{l}
\psi_i^T= \varphi^T_i \mbox{ for } i=1,\ldots,N+1 \\
\psi_{N+2}^T = [{\bf x}_T-{\bf x}]^{\mbox{t}}{\mathcal K}^{-1}{\bf w}_h + 1.
\end{array}
\right.
\end{equation}
\indent Next we define, \\

\noindent \underline{Hermite finite element spaces $U_h$ and $V_h$:}\\
\noindent Consider that for every interface $F$ (an inner edge for $N=2$ and an inner face for $N=3$) of two 
elements in ${\mathcal T}_h$,  ${\bf n}_F$ is oriented in the same manner for both of them.  
Then every function in $v \in V_h$ (resp. $\in U_h$) is such that its restriction to every $T \in {\mathcal T}_h$ is 
a $N+2$ coefficient quadratic function of the form (\ref{v}),  
whose degrees of freedom of the type ${\mathcal V}_F$ (resp. ${\mathcal U}_F)$ coincide on both sides of every interface $F$ 
of a pair of elements in ${\cal T}_h$. \\

We proceed by setting the discrete variational problem (\ref{pbh}) below, aimed at approximating (\ref{model}), 
whose bi-linear form $a_h$ and linear form $L_h$ are given by (\ref{ahLps}):

Find $ u_h \in U_h$ such that  for all $v \in V_h $ \\
\begin{equation}
\label{pbh}
a_h(u_h,v) = L_h(v),
\end{equation}

\noindent holds, where $\forall u \in U_h$ and $\forall v \in V_h$,
    
\begin{equation}
\label{ahLps}
\left\{
\begin{array}{ll}
a_h(u,v):= & \displaystyle \sum_{T \in {\cal T}_h} 
\displaystyle [ (\nabla \cdot {\cal K} \nabla u - {\bf w}_h^1 \cdot \nabla u, \Pi_T[v])_T \\ [ 2ex]
& +(\nabla u,{\cal K} \nabla v+{\bf w}_h \Pi_T[v])_{T} + (u,\nabla \cdot {\cal K} \nabla v)_T ];\\ [ 2ex]
L_h(v):= & -(f,\Pi_h[v]).
\end{array}
\right.
\end{equation}
Now let us  consider the space 
\[ V := \; \{v|v \in H^1(\Omega);\nabla \cdot {\cal K} \nabla v \in L^2(\Omega)\}. \]  
Clearly $a_h$ can be extended to $(U_h+V)\times(V_h+V)$. Then  
we further introduce the functional $\parallel \cdot \parallel_h: U_h+V_h+V \longrightarrow \Re$ given by:
\begin{equation}
\parallel v \parallel_{h} := \; \displaystyle \left[(\Pi_h v,\Pi_h v) + \sum_{T \in {\cal T}_h} 
\left\{ (\nabla v, \nabla v)_{T} +  
(\nabla \cdot {\cal K} \nabla v , \nabla \cdot {\cal K} \nabla v)_T \right\} \right]^{1/2}. 
\end{equation}
The expression $\parallel \cdot \parallel_h$ obviously defines a norm over $V$, $U_h$ and $V_h$. 
In this manner, it is not difficult to establish the continuity of $a_h$ over $(U_h + V) \times (V_h + V)$ 
with a mesh independent constant $M$ (cf. the proof of \textbf{Proposition 3.1} hereafter):
\begin{equation}\label{contah}
a_h(u, v) \le M \parallel  u \parallel_h  \,  \parallel  v \parallel_h.
\end{equation}
On the other hand there is no way for $a_h$ to be coercive. Hence we resort to an inf-sup condition for $a_h$ over 
$U_h \times V_h$ \cite{Babuska}, which directly implies that (\ref{pbh}) has a unique solution. 
More specifically the following stability result was proved in \cite{ICNAAM2013}. 
\begin{e-proposition} \label{stabil}
(\cite{ICNAAM2013}) If $h$ is sufficiently small and ${\bf w} \in [W^{1,\infty}(\Omega)]^N$, 
there exists a constant $\alpha>0$ independent of $h$ such that 
\begin{equation} 
\label{infsup}
 \forall u \in U_h \; \displaystyle \sup_{v \in V_h \setminus \{0\}} \frac{a_h(u,v)}{\parallel v \parallel_h} 
 \geq \alpha 
 \parallel u \parallel_h.
\end{equation}
\end{e-proposition}

In the next section we derive estimates for $\parallel u - u_h \parallel_h$ using a modified Strang Lemma for non coercive 
problems given in \cite{COAM}. In this aim we have to consider the following auxiliary problem:

Find $ u_h^{*} \in U_h$ such that for all $v \in V_h $ \\
\begin{equation}
\label{pbhstar}
a^{*}_h(u^{*}_h,v) = L_h(v),
\end{equation}

\noindent holds, where $\forall u \in U_h + V$ and $\forall v \in V_h + V$,
    
\begin{equation}
\label{ahstar}
\left\{
\begin{array}{ll}
a^{*}_h(u,v):= & \displaystyle \sum_{T \in {\cal T}_h} 
\displaystyle [ (\nabla \cdot {\cal K} \nabla u - {\bf w} \cdot \nabla u, \Pi_T[v])_T \\ [ 2ex]
& 
+(\nabla u,{\cal K} \nabla v+{\bf w}_h \Pi_T[v])_{T} + (u,\nabla \cdot {\cal K} \nabla v)_T ].
\end{array}
\right.
\end{equation} 

Similarly to the case of problem (\ref{pbh}) (cf. \cite{ICNAAM2013}) we can prove,  
\begin{theorem}
\label{ExUn}
Problem (\ref{pbhstar}) has a unique solution and moreover there exists a constant $C^{*}$ independent of $h$ such that 
\begin{equation}
\label{bounduhstar}
\parallel u_h^{*} \parallel_{h} \leq C^{*} \parallel f \parallel_{0,\Omega}.
\end{equation} 
\end{theorem}

Before proving \textbf{Theorem 2.2} we establish a stability result for problem (\ref{pbhstar}), namely,

\begin{e-proposition}
\label{stabilstar}
If $h$ is sufficiently small and ${\bf w} \in [W^{1,\infty}(\Omega)]^N$, 
there exists a constant $\alpha^{*}>0$ independent of $h$ such that 
\begin{equation} 
\label{infsupstar}
 \forall u \in U_h \; \displaystyle \sup_{v \in V_h \setminus \{0\}} \frac{a_h^{*}(u,v)}{\parallel v \parallel_h} 
 \geq \alpha^{*} \parallel u \parallel_h.
\end{equation}
\end{e-proposition}
	
\noindent \underline{Proof:} 
Given $u \in U_h$ define $v:=v_1+v_2+v_3$, where $v_i \in V_h$ for $i=1,2,3$ are defined as follows:\\ 
$v_1= \theta_1 w_1$, $\theta_1$ being a non negative constant to be 
specified, and $w_1$ being defined by $\Pi_h[w_1] = \Pi_h[u]$, together with 
$({\mathcal K} \nabla w_1 + {\bf w}_h \Pi_h[w_1]) \cdot {\bf n}_T = ({\mathcal K} \nabla u) \cdot {\bf n}_T$ for every $T 
\in {\mathcal T}_h$, where ${\bf n}_T$ is the outer normal on $\partial T$. \\
$v_2$ equals $\theta_2 \; \nabla \cdot {\cal K} \nabla u$ in every $T \in {\cal T}_h$, where $\theta_2$ 
is a non negative constant to be specified. \\
$v_3$ is constructed by applying Theorem 4 of \cite{RaviartThomas}. According to it there exists a field 
${\bf p} \in {\bf Q}_h:= \{{\bf q}\; | \; \exists u \in U_h  
\mbox{ such that } {\bf q}_{|T} = {\cal K} \nabla 
u_{|T} \; \forall T \in {\cal T}_h \}$, satisfying for a constant $\tilde{C}$ independent of $h$:

\begin{equation}
\label{intermed2}
\nabla \cdot {\bf p} = \Pi_h[u] \mbox{ in } \Omega; \hspace{1cm}
\parallel {\bf p} \parallel_{0,\Omega} \leq \tilde{C} 
\parallel \Pi_h[u] \parallel_{0,\Omega}.
\end{equation}

\noindent Then recalling that the normal traces over the faces of the elements in ${\cal T}_h$ of fields belonging 
to ${\bf Q}_h$ are constant \cite{RaviartThomas} $v_3$ is defined in such a way that  $\forall T \in {\mathcal T}_h$, 
$({\mathcal K} \nabla  v_3 + {\bf w}_h \Pi_T[v_3]) \cdot {\bf n}_T={\bf p} \cdot {\bf n}_T$ $\forall T \in {\cal T}_h$ and  
$\Pi_h[v_3] = - \theta_1 \Pi_h[u]$.\\
\indent It is clear that $\nabla \cdot {\mathcal K} \nabla w_1 = \nabla \cdot {\mathcal K} \nabla u$. Moreover 
by construction we have, 
\begin{equation}
 \label{intermed2bis}
 \oint_{\partial T} {\mathcal K} \nabla    w_1 \cdot {\bf n}_T w_1 \;dS - (\nabla \cdot {\mathcal K} w_1,w_1)_{T} = 
 ({\mathcal K}\nabla    w_1,\nabla    w_1)_{T} 
 = ({\mathcal K}\nabla    u-{\bf w}_h \Pi_T[u],\nabla w_1)_{T}  \; \forall T \in {\mathcal T}_h.
 \end{equation}
 
 \noindent Then $\lambda$ and $\Lambda$ being the smallest and the largest eigenvalue of ${\cal K}$, 
 after straightforward manipulations it follows that 
\begin{equation}
\label{intermediary}
\parallel \nabla    w_1 \parallel_{0,T} \leq \lambda^{-1} (\Lambda \parallel \nabla    u \parallel_{0,T} + 
\parallel {\bf w} \parallel_{0,\infty,\Omega} \parallel \Pi_T[u] \parallel_{0,T} ).
\end{equation}

\noindent This implies that for $\tilde{C}_1=\lambda^{-1}(\Lambda^2+\parallel {\bf w} \parallel_{0,\infty,\Omega}^2)^{1/2}$, 
$\displaystyle \sum_{T \in {\mathcal T}_h} \parallel \nabla    w_1 \parallel_{0,T}^2 \leq  
\tilde{C}_1^2 \parallel u \parallel_h^2$, which immediately yields,
\begin{equation}
 \label{intermediarybis}
\parallel v_1 \parallel_h \leq C_1 \parallel u \parallel_h, \mbox{ with } C_1 = \theta_1 (1+\tilde{C}_1^2)^{1/2}.
\end{equation}

As for $v_2$ we readily have, 
\begin{equation}
 \label{intermediaryter}
\parallel v_2 \parallel_h \leq C_2 \parallel u \parallel_h, \mbox{ with } C_2 = \theta_2.
\end{equation}

On the other hand by construction $v_3$ fulfills $\nabla \cdot [{\mathcal K} \nabla v_3]_{|T} = \nabla \cdot {\bf p}_{|T}=\Pi_T[u]$, 
$\forall T \in {\mathcal T}_h$, and hence,    
\begin{equation}
 \label{intermed2ter}
\lambda \parallel \nabla    v_3 \parallel_{0,T}^2 \leq ({\mathcal K} \nabla    v_3,\nabla    v_3)_T 
=\oint_{\partial T} ({\mathcal K} \nabla v_3) \cdot {\bf n}_T v_3  dS -(v_3,\nabla \cdot {\bf p})_{T} 
= (\theta_1 {\bf w}_h \Pi_T[u]+{\bf p},\nabla v_3)_{T}
\end{equation}

\noindent It easily follows that
\begin{equation}
\label{intermed3}
\parallel v_3 \parallel_h  \leq C_3 \parallel \Pi_h[u] \parallel_{0,\Omega} \leq C_3 \parallel u \parallel_h, \mbox{ where } 
C_3=[(\tilde{C}+ \theta_1 \parallel {\bf w} \parallel_{0,\infty,\Omega})^2\lambda^{-2}+\theta_1^2+1]^{1/2}. 
\end{equation} 

Now taking into account (\ref{intermed2}) and (\ref{intermed2bis}), after straightforward calculations we obtain,
\begin{equation}
\label{intermed1}
a_h^{*}(u,v_1) = \theta_1 \displaystyle \sum_{T \in {\cal T}_h}
\{ 2( \nabla \cdot {\cal K} \nabla u,\Pi_T[u])_{T} + ({\mathcal K} \nabla u,\nabla u)_{T}  
-({\bf w} \cdot \nabla    u, \Pi_T[u])_T \};
\end{equation}
\begin{equation}
\label{intermed1bis}
a_h^{*}(u,v_2) = \theta_2 \displaystyle \sum_{T \in {\cal T}_h} 
\{\parallel \nabla \cdot {\cal K} \nabla u \parallel_{0,T}^2 + 
([{\bf w}_h-{\bf w}] \cdot \nabla  u,\nabla \cdot {\cal K} \nabla u)_T \};
\end{equation}
\begin{equation}
\label{intermed1ter}
a_h^{*}(u,v_3) =  \displaystyle \sum_{T \in {\cal T}_h} \{ \parallel \Pi_T[u] \parallel_{0,T}^2 + 
\theta_1 [  ({\bf w} \cdot \nabla    u,\Pi_T[u] )_T -(\nabla \cdot {\cal K} \nabla u,\Pi_T[u])_T ] 
+({\bf p},\nabla u)_T \}.
\end{equation}


\noindent Then, recalling the definition of ${\bf w}_h$, there exists a mesh independent constant $C_W$ (cf. \cite{Ciarlet}) such that  
$\parallel {\bf w} - {\bf w}_h \parallel_{0,\infty,T} \leq C_W h | {\bf w} |_{1,\infty,T}$ $\forall T \in {\mathcal T}_h$. Using this fact, together with (\ref{intermed1})-(\ref{intermed1bis})-(\ref{intermed1ter}), simple manipulations lead to:
\begin{equation}
\label{intermed4}
\left\{
\begin{array}{l}
a_h^{*}(u,v) \geq \parallel \Pi_h[u] \parallel_{0,\Omega}^2 + \displaystyle \sum_{T \in {\cal T}_h} 
\displaystyle \left[ \theta_1 \lambda \parallel \nabla   u \parallel_{0,T}^2 + \frac{\theta_2}{2} 
\parallel \nabla \cdot {\cal K} \nabla u \parallel_{0,T}^2 - \parallel {\bf p} \parallel_{0,T} \parallel \nabla u \parallel_{0,T} \right] \\
- \displaystyle \sum_{T \in {\cal T}_h} \left[ \theta_1 \parallel \nabla \cdot {\cal K} \nabla u \parallel_{0,T} \parallel \Pi_T[u] \parallel_{0,T}) + \theta_2 C_W^2 h^2 | {\bf w} |_{1,\infty,\Omega}^2 \parallel \nabla   u \parallel_{0,T}^2/2 \right] \\
\geq \parallel \Pi_h[u] \parallel_{0,\Omega}^2/4 + \displaystyle \sum_{T \in {\cal T}_h} 
\{(\theta_1 \lambda - \tilde{C}^2 - \theta_2 C_W^2 h^2 |{\bf w}|_{1,\infty,\Omega}^2/2)
\parallel \nabla   u \parallel_{0,T}^2 \\ + (\theta_2/2-\theta_1^2/2)\parallel \nabla \cdot {\cal K} \nabla u \parallel_{0,T}^2 \}
\end{array}
\right.
\end{equation} 

\noindent Now if we assume that $h^2 \leq \beta (C_W |{\bf w}|_{1,\infty,\Omega})^{-2}$ with $\beta \leq 
4 \lambda^2/[D+(D^2+8 \lambda^2)^{1/2}]$ for $D=1+4\tilde{C}^2$, we may choose $\theta_1>0$ satisfying $\theta_1 \lambda -\tilde{C}^2 - \beta \theta_2/2 
\geq 1/4$ with $\theta_2 =1/2+\theta_1^2$. It follows from (\ref{intermed4}),   
(\ref{intermediarybis}), (\ref{intermediaryter}) and (\ref{intermed3}), that,
\begin{equation}
\label{infsup2}
a_h^{*}(u,v) \geq  \parallel u \parallel_{h}^2/4; \hspace{1cm} 
\parallel v \parallel_h \leq C \parallel u \parallel_h, \mbox{ with } C=[3(C_1^2+C_2^2+C_3^2)]^{1/2}.   
\end{equation}
This immediately yields (\ref{infsupstar}) with $\alpha^{*}=1/(4C)$. \rule{2mm}{2mm}\\

\noindent \underline{Proof of \textbf{Theorem 2.2}:} Since $V_h$ is a finite dimensional space, according to \cite{Babuska} the existence and uniqueness of a solution to (\ref{pbhstar}) follows from (\ref{infsupstar}). 
Moreover, combining (\ref{pbhstar}) and (\ref{infsupstar}) we easily obtain, 
\begin{equation} 
\label{boundf}
\alpha^{*} \parallel u_h^{*} \parallel_{h} \leq \displaystyle \sup_{v \in V_h \setminus \{0\}} \frac{L_h(v)}{\parallel 
\Pi_h[v] \parallel_{0,\Omega}}.
\end{equation}
Since $L_h(v) = \int_{\Omega} f \Pi_h[v] dx$ from (\ref{boundf}) we finally derive (\ref{bounduhstar}) with $C^{*}=[\alpha^{*}]^{-1}$. \rule{2mm}{2mm} \\

\section{Convergence results}

Henceforth we denote by $\nabla_h$ the operator from $V+U_h+V_h$ onto $L^2(\Omega)$ defined by 
\[ [\nabla_h w]_{|T}=\nabla [w_{|T}] \; \forall T \in {\mathcal T}_h, \; \forall w \in V+U_h+V_h.\]
Notice that for any function $u \in V+U_h$, $\nabla \cdot {\mathcal K} \nabla u$ is well-defined in $L^2(\Omega)$ (cf. \cite{Thomas}) and hence there is no need to use the operator $\nabla_h$ in this case. \\   

In order to study the convergence of $u_h$ to $u$ in appropriate norms we first note that 
from the properties of $V_h$ and equation (\ref{model}) we easily infer that $u$ satisfies 

\begin{equation} 	
\label{ahstaru}
a^{*}_h(u,v) = L_h(v) \; \forall v \in V_h
\end{equation}
 
From the continuity of $a^{*}_h$ and the uniform stability result proved in \textbf{Proposition 2.3}, we may apply the generalized First and Second Strang's inequality for the weakly coercive case, namely, inequality (32) of \cite{COAM}. 
In the case under study this writes,

\begin{equation}
\label{StrangCOAM}
\parallel u - u_h \parallel_h \leq \displaystyle \frac{M^{*}}{\alpha^{*}} \displaystyle \inf_{w \in U_h} \parallel u -w \parallel_h 
+ \displaystyle \frac{1}{\alpha} \displaystyle \sup_{v \in V_h \setminus \{ 0 \} } \frac {[a_h^{*}-a_h](u^{*}_h,v)}{\parallel v \parallel_h}
\end{equation}
\noindent where $M^{*}$ is a constant such that 
\begin{equation} 
\label{Mstar}
a^{*}_h(u,v) \leq M^{*} \parallel u \parallel_h \parallel v \parallel_h \; \forall u \in V + U_h, \; \forall v \in V+V_h. 
\end{equation}
\begin{e-proposition}
There exists a constant $M^{*}$ independent of $h$ such that (\ref{Mstar}) holds.
\end{e-proposition}

\noindent \underline{Proof:} Since $(u,\nabla \cdot {\mathcal K} \nabla v)_T = (\Pi_T[u], \nabla \cdot {\mathcal K} \nabla v)_T
\; \forall T \in {\mathcal T}_h$ and $\forall v \in V_h$, 
we trivially have,
\begin{equation}
\label{Mstar1}
\left\{
\begin{array}{l}
a^{*}_h(u,v) \leq (\parallel \nabla \cdot {\mathcal K} \nabla_h u \parallel_{0,\Omega}  + 
\parallel {\bf w} \parallel_{0,\infty,\Omega} \parallel \nabla_h u \parallel_{0,\Omega}) \parallel \Pi_h[v] \parallel_{0,\Omega} + \\
\parallel \nabla_h u \parallel_{0,\Omega} (\Lambda \parallel \nabla_h v \parallel_{0,\Omega} + \parallel {\bf w} \parallel_{0,\infty, \Omega} \parallel \Pi_h[v] \parallel_{0,\Omega}) + \displaystyle \sum_{T \in {\mathcal T}_h} 
\parallel \Pi_T[u] \parallel_{0,T} \parallel \nabla \cdot {\mathcal K} \nabla v \parallel_{0,T}. 
\end{array}
\right.
\end{equation}
(\ref{Mstar1}) immediately yields (\ref{Mstar}) with $M^{*} = 2 \max[1,2 \parallel {\bf w} \parallel_{0,\infty,\Omega}, \Lambda]$.
\rule{2mm}{2mm}\\

Next we prove the validity of the following a priori error estimate for the method under study:

\begin{theorem}
\label{theorem1} 
Assume that ${\bf w} \in [W^{1,\infty}(\Omega)]^N$ and   
$h$ is sufficiently small. Then if $u \in H^2(\Omega)$ and $f \in H^1(\Omega)$ there exists a mesh independent constant $C^{'}$ 
such that, 
\begin{equation}
\label{converge1}
\parallel u - u_h \parallel_h \leq C^{'} h [ \parallel u \parallel_{2,\Omega} + \parallel f \parallel_{1,\Omega}]
\end{equation}
\end{theorem}

\noindent \underline{Proof:} By standard results applying to the $RT_0$ method, and since 
$\parallel \Pi_h[u - u_h] \parallel_{0,\Omega}$ is obviously bounded above by a mesh independent constant times $h|u|_{1,\Omega}$, 
for a suitable constant $C_I$ independent of $h$ it holds,  
\begin{equation}
\label{errorestim1}
\displaystyle \inf_{w \in V_h} \parallel u - w \parallel_h \leq C_I h [ \parallel u \parallel_{2,\Omega} + |f|_{1,\Omega}]
\end{equation} 
On the other hand we have $|[a^{*}_h-a_h](u^{*}_h,v)|=|([{\bf w}-{\bf w}^1_h] \cdot \nabla_h u^{*}_h,\Pi_h[v])|$. Hence 
for a mesh independent constant $C^{*}_W$ such that $\parallel {\bf w}-{\bf w}^1_h \parallel_{0,\infty,\Omega} \leq C^{*}_W h 
| {\bf w} |_{1,\infty,\Omega}$ we derive,

\begin{equation} 
\label{errorestim2}
|[a^{*}_h-a_h](u^{*}_h,v)| \leq C^{*}_W h |{\bf w}|_{1,\infty,\Omega} 
\parallel \nabla_h u^{*}_h \parallel_{0,\Omega} \parallel v \parallel_h.
\end{equation}
Taking into account (\ref{bounduhstar}), (\ref{StrangCOAM})-(\ref{errorestim1})-(\ref{errorestim2}) readily yield 
(\ref{converge1}), $C^{'}$ being a mesh independent constant. \rule{2mm}{2mm} \\ 

Next we give a fundamental result of this work:

\begin{theorem}
\label{theorem2}
If $\Omega$ is convex, ${\bf w} \in [W^{2,4}(\Omega)]^N$ and $h$ is sufficiently small, there exists a constant $C^{''}$ independent of $h$ such that, 
\begin{equation}
\label{converge2} 
\parallel u - u_h \parallel_{0,\Omega} \leq C^{''} h^2 \; ( \parallel u \parallel_{2,\Omega} + \parallel f \parallel_{1,\Omega}).  
\end{equation}
\end{theorem}

\noindent \underline{Proof:}  
The proof is a non-trivial extension of the proof of Theorem 2.3 in \cite{CAM2013}, where the quadratic convergence was shown for 
the pure diffusion case. We first observe that $\bw \in [W^{1,\infty}(\Omega)]^N$ according to the Sobolev Embedding Theorem 
\cite{Adams}.  Moreover using the definitions of $a_h$ and $\Pi_h$, together with the continuity of the normal components of 
${\mathcal K} \nabla v_h+{\bf w}_h \Pi_h [v_h]$ on $\partial T$ for $v_h \in V_h$, we easily obtain,
\begin{equation}\label{proof_eq:1}
a_h(u-u_h, v_h) + ((\bw^1_h - \bw ) \cdot \nabla u, \Pi_h[v_h]) = 0,
\end{equation}
for all $v_h \in V_h$. Similarly, owing to the continuity of the normal components of ${\mathcal K} \nabla u_h$ on $\partial T$ 
(cf. \cite{CAM2013}):
\begin{equation}\label{proof_eq:2}
(u-u_h,  \nabla \cdot {\cal K} \nabla v) = a_h(u-u_h, v) + (\nabla 
\cdot {\mathcal K} \nabla (u-u_h), v-\Pi_h[v])+((\bw^1_h-\bw_h) \cdot \nabla_h (u - u_h), \Pi_h[v]), 
\end{equation}
for all $v \in \{v  |  v \in H^1_0(\Omega), \nabla \cdot {\cal K} \nabla v \in L^2(\Omega) \}$. By using the Aubin-Nitsche trick 
and \eqref{proof_eq:2} we can write 
\begin{equation}
\begin{array}{l}
\parallel u - u_h \parallel_{0, \Omega} \; \; = \displaystyle \sup_{v \in D(\Omega) \backslash \{0\} }  
\dfrac{(u -u_h, \nabla \cdot {\cal K} \nabla v)}{\parallel \nabla \cdot {\cal K } \nabla v \parallel_{0, \Omega}} \;  = \\
\displaystyle \sup_{v \in D(\Omega) \backslash \{0\} } \dfrac{ a_h(u-u_h, v) + 
(\nabla \cdot {\mathcal K} \nabla (u-u_h), v-\Pi_h[v]) + ((\bw^1_h  -  \bw_h)\cdot \nabla_h (u - u_h), \Pi_h[v]) }{\parallel 
\nabla \cdot {\cal K} \nabla v \parallel_{0, \Omega}}
 \label{proof_eq:3}
\end{array}
\end{equation}
where $D(\Omega):= V \cap H^1_0(\Omega). $ Let now $u \neq u_h$ (if $u=u_h$ then \eqref{converge2} trivially holds) and 
$B_D({\bf 0},1):=\{v  |  v \in D(\Omega), \; \parallel \nabla \cdot {\cal K} \nabla v \parallel_{0,\Omega}=1\}$. 
We know that there exists $v_0 \in B_D({\bf 0},1)$ such that $\nabla \cdot {\cal K} \nabla v_0\label{proof_eq:8bis} = 
\dfrac{u-u_h}{ \parallel u-u_h \parallel_{0,\Omega}}$ (cf. \cite{Evans}). 
Thus it is easy to see that due to \eqref{proof_eq:3} it holds,
\begin{equation}\label{proof_eq:4}
\parallel u - u_h \parallel_{0, \Omega} =  a_h(u-u_h, v_0) + (\nabla \cdot {\mathcal K} \nabla (u-u_h), v_0-\Pi_h[v_0])+
((\bw^1_h - \bw_h )\cdot \nabla_h (u - u_h), \Pi_h[v_0]).
\end{equation}
By combining \eqref{proof_eq:1} and \eqref{proof_eq:4} we further get for any $v_h \in V_h$
\begin{eqnarray}\label{proof_eq:5}
\parallel u - u_h \parallel_{0, \Omega} =  a_h(u-u_h, v_0-v_h) +(\nabla \cdot {\mathcal K} \nabla (u-u_h), v_0 -\Pi_h[v_0]) \\
[ 2ex] \nonumber
+ ((\bw^1_h - \bw_h )\cdot \nabla_h (u - u_h), \Pi_h[v_0])-((\bw^1_h - \bw)\cdot \nabla u, \Pi_h[v_h]).
\end{eqnarray}
Since $D(\Omega) \subset H^2(\Omega)$ in case $\Omega$ is convex (cf. \cite{Evans}), 
we can define the standard interpolate $I_h v \in V_h$ of every $v \in D(\Omega)$, based on the degrees of freedom of $V_h$. 
Taking $v_h = I_h v_0$ in \eqref{proof_eq:5}, we get,
\begin{eqnarray}\label{proof_eq:6}
\parallel u - u_h \parallel_{0, \Omega} =  a_h(u-u_h, v_0-I_h v_0) +(\nabla \cdot {\mathcal K} \nabla (u-u_h), v_0-\Pi_h[v_0])\\
[ 2ex] \nonumber
+ ((\bw^1_h - \bw_h )\cdot \nabla_h (u - u_h), \Pi_h[v_0])-((\bw^1_h - \bw)\cdot \nabla u, \Pi_h[I_h v_0]).
\end{eqnarray}
By using the continuity \eqref{contah} of $a_h$ and the definition of $v_0$, together with the approximation properties of $I_h$ for $h$ sufficiently small 
(see \cite{CAM2013}, p. 239 for details), we have for a suitable $h$-independent constant $\underline{C}$,
\begin{equation}\label{proof_eq:7}
a_h(u-u_h, v_0-I_h v_0) \leq \displaystyle \frac{3}{4} \parallel u - u_h \parallel_{0,\Omega} + 
\underline{C} \parallel u - u_h \parallel_h [\parallel v_0 - I_h v_0 \parallel_{0,\Omega}+\parallel \nabla_h(v_0 - I_h v_0) \parallel_{0,\Omega}].
\end{equation}  
\noindent Applying to \eqref{proof_eq:6} and \eqref{proof_eq:7} standard results to estimate $\parallel v_0 - 
\Pi_h[v_0] \parallel_{0,\Omega}$ together with $\parallel v_0 - I_h v_0 \parallel_{0,\Omega}+ 
\parallel \nabla_h(v_0 - I_h v_0) \parallel_{0,\Omega}$, recalling \eqref{converge1} we easily 
conclude  that there exists another constant $\bar{C}$ independent of $h$ such that, 
\begin{equation}\label{proof_eq:8}
\parallel u - u_h \parallel_{0, \Omega} \leq  \bar{C}  h^2 ( | u |_{2,\Omega} + | f |_{1,\Omega}) + 4[|((\bw^1_h - 
\bw_h )\cdot \nabla_h 
(u - u_h), \Pi_h[v_0])|+|((\bw^1_h - \bw)\cdot \nabla u, \Pi_h[I_h v_0])|].
\end{equation}
We proceed by estimating the two terms in brackets on the right hand side of \eqref{proof_eq:8} denoted by $T_1$ and $T_2$. \\
First we note that from the Sobolev Embedding Theorem and the convexity of $\Omega$ (cf. \cite{Grisvard}), there exist 
constants $C_{\infty}$ and $C_{\infty}^{'}$ depending only on $\Omega$ such that 
$$\parallel \Pi_h[v_0] \parallel_{0,\infty,\Omega} \leq \parallel v_0 \parallel_{0,\infty,\Omega} \leq C_{\infty} 
\parallel v_0 \parallel_{2,\Omega} \leq C_{\infty}^{'}.$$
Thus for a mesh independent constant $C_2$ we have:
\begin{equation}\label{proof_eq:9}
T_1 \le \parallel \bw^1_h - \bw_h  \parallel_{0,\Omega} \,  \parallel \nabla_h(u - u_h) \parallel_{0,\Omega} \,  
\parallel \Pi_h[v_0]  \parallel_{0,\infty,\Omega} \leq C_2 h^2 | u |_{2,\Omega}.
\end{equation} 
For deriving the estimate (\ref{proof_eq:9}) we used the fact that both $\bw^1_h$ and $\bw_h$ are interpolates of $\bw$, the result  
\eqref{converge1} and the boundedness of  $\parallel \Pi_h[v_0] \parallel_{0,\infty,\Omega}$. 
$C_2$ is the product of $C_{\infty}^{'}$ with another constant not depending on $h$ and the semi-norm of 
$\bf w$ in $[H^{1}(\Omega)]^N$. \\
The term $T_2$ can be estimated in a similar way. Using now the fact that $\bw^1_h$ is a piecewise linear 
interpolate of $\bw$, the regularity of the solution $u$ and standard properties of $\Pi_h$ and $I_h$, there holds

\begin{equation}
\label{proof_eq:10}
T_2 \le \parallel \bw^1_h - \bw  \parallel_{0,4,\Omega} \,  \parallel \nabla u \parallel_{0,4,\Omega} \,  
\parallel \Pi_h[I_h v_0]  \parallel_{0,\Omega} \le C_3 h^2 \parallel u \parallel_{2,\Omega}, 
\end{equation}
where $C_3$ equals an $h$-independent constant times $| \bw |_{2,4,\Omega}$. Putting together 
\eqref{proof_eq:8}-\eqref{proof_eq:9}-\eqref{proof_eq:10}, we obtain the quadratic convergence 
\eqref{converge2}. \rule{2mm}{2mm}\\ 

\section{A variant for the equations in divergence form}

\hspace{4mm} Like in \cite{Douglas} it is possible to consider a variant of the method described in Section 2 applying to the case 
where the normal component of the total flux $-{\mathcal K} \nabla u + {\bf w} u$ is continuous across the element interfaces. 
In the case of the mixed formulation this corresponds to introducing the auxiliary variable ${\bf p}$ given by the above expression, 
and write the C-D equation equation (\ref{model}) in divergence form, namely

Find $ u $ satisfying $ u=0 $  on $\Gamma$ and $ {\bf p} $ such that \\
\begin{equation}
\label{model2}
\left\{
\begin{array}{cccc}
\nabla \cdot {\bf p} - \nabla \cdot {\bf w} \; u &= &f & \mbox{ in } \Omega, \\[ 2ex]
{\bf p} +{\mathcal K} \nabla u - {\bf w} u & = & 0 & \mbox{ in } \Omega.
\end{array}
\right.
\end{equation}

\noindent Recalling the space ${\bf H}(div;\Omega):= \{{\bf q}\; | \; {\bf q} \in [L^2(\Omega)]^N, \; 
\nabla \cdot {\bf q} \in L^2(\Omega) \}$, a natural weak (variational) formulation equivalent to system (\ref{model2}) is given 
in \cite{Douglas}, that is,

Find  $u \in L^2(\Omega) $ and ${\bf p} \in {\bf H}(div;\Omega)$ such that for all $v \in  L^2(\Omega)$, and for all 
$ {\bf q} \in {\bf H}(div;\Omega)$, \\
\begin{equation}
\label{model3}
\left\{
\begin{array}{rcl}
(\nabla \cdot {\bf p},v) - (\nabla \cdot {\bf w} \; u,v)  &=& (f,v) \\[ 2ex]
({\mathcal K}^{-1} {\bf p}, {\bf q}) 
- (u, \nabla \cdot {\bf q}) - ({\mathcal K}^{-1} {\bf w} \; u, {\bf q}) &=& 0.
\end{array}
\right.
\end{equation} 

The extension of $RT_0$ to the C-D equation considered in \cite{Douglas} consists of using the Raviart-Thomas interpolation of 
the lowest order to represent ${\bf p}$ and ${\bf q}$ - i.e. to approximate ${\bf H}(div;\Omega)$ -, and the space of 
constant functions in each element of the partition ${\mathcal T}_h$ to represent $u$ and $v$. In contrast, here we shall mimic 
(\ref{model3}) by resorting to the space $U_h$, after adding up both relations in (\ref{model3}). More specifically  
we take in each element $T \in {\mathcal T}_h$, ${\bf q}_{|T}={\mathcal K} \nabla v_{|T}$ for $v \in U_h$. Now $u_h$ will be searched 
for in a space $W_h$ defined hereafter. First we have to    
construct field $\tilde{\bf w}_h$ to replace ${\bf w}_h$ (cf. Section 2), in order to preserve optimality of the approximation of 
$u$. In this aim it suffices that $\tilde{\bf w}_h$ be of the form $c{\bf x} + {\bf d}$ in each 
$T \in {\mathcal T}_h$ for suitable real number $c$ and real vector ${\bf d}=[d_1,\ldots,d_N]^{{\mbox t}}$. This representation is compatible with the requirement that the normal component of the flux variable ${\bf p}= -{\mathcal K} \nabla + 
{\bf w} u$ be continuous across the mesh edges at discrete level. The natural choice of $\tilde{\bf w}_h$ is 
certainly the interpolate of ${\bf w}$ in the Raviart-Thomas ($RT_0$) space. Now we define the\\

\noindent \underline{Hermite finite element space $W_h$:}\\
$W_h$ is the space of functions $v$ of the form 
${\bf x}^{\mbox{t}}\{{\mathcal K}^{-1}[a{\bf x}/2 + {\bf b}]\}+d$ in every $T \in {\mathcal T}_h$, such that the mean normal flux  
$\int_F (-{\mathcal K} \nabla v + \tilde{\bf w}_h \Pi_h[v])\cdot {\bf n}_{F} ds/meas(F)$ is continuous across all the inner edges 
or faces $F$ of the partition. This is about all that is needed to complete the definition of $W_h$. Indeed using a procedure 
very similar to the one in Section 2 (cf. (\ref{basisU})) it is possible to uniquely determine the $N+2$ local basis functions 
$\eta^T_i$ for each $N$-simplex $T \in {\mathcal T}_h$ related to the above set of degrees of freedom completed 
with the function mean value in $T$, defining $W_h$ locally. More precisely, 
setting $[\tilde{\bf w}_h]_{|T} = a_w^T {\bf x} + {\bf b}_w^T$, since by definition 
$\Pi_T \eta_i^T=0$ for $i=1,\ldots,N+1$ and $\Pi_T \eta^T_{N+2}=1$, recalling (\ref{basisU}) we have: 
\begin{equation}
\label{basisW}
\left\{
\begin{array}{l}
\eta_i^T = \varphi^T_i \mbox{ for } i=1,\ldots,N+1, \\
\eta^T_{N+2}={\bf x}^{\mbox{t}}\{{\cal K}^{-1}[a_w^T {\bf x}/2 + {\bf b}_w^T]\} +  
1-\int_T {\bf x}^{\mbox{t}}\{{\cal K}^{-1}[a_w^T {\bf x}/2 + {\bf b}_w^T]\}dx/meas(T).
\end{array}
\right.
\end{equation}
\indent Now we replace in (\ref{model3}) :

\begin{itemize}
\item 
$u$ with $\Pi_h[u_h]$;
\item 
${\bf w}$ with $\tilde{\bf w}_h$;
\item
${\bf p}$ with $-{\mathcal K} \nabla_h u_h +\tilde{\bf w}_h \Pi_h[u_h]$ (taking $u_h \in W_h$); 
\item
${\bf q}$ with $-{\mathcal K} \nabla_h v$ (taking $v \in U_h$);
\item
$f$ with $\Pi_h[f]$.
\end{itemize}

This leads to the following equation:

\begin{equation}
 \label{eqn1}
\begin{array}{l}
  \displaystyle \sum_{T \in {\mathcal T}_h} \left[ (\nabla \cdot \{{\mathcal K}\nabla u_h -\tilde{\bf w}_h 
  \Pi_T[u_h]\},v)_{T} + (\nabla \cdot \tilde{\bf w}_h \Pi_T[u_h],v)_{T} + \right. \\
 ({\mathcal K} \nabla u_h-\tilde{\bf w}_h \Pi_T[u_h], \nabla v)_T + (\Pi_T[u_h], \nabla \cdot {\mathcal K} \nabla v)_{T} + 
 \\[ 2ex]
 \left. (\tilde{\bf w}_h \Pi_T[u_h], \nabla v)_{T} \right] = -(\Pi_h[f],v) \; \forall v \in  U_h.
\end{array}
\end{equation}

After straightforward simplifications, and taking into account that $(\Pi_h[f],v)=(f,\Pi_h[v])$, 
we come up with the following Hermite finite element counterpart of (\ref{model}):

Find $u_h \in W_h$ such that for all  $v \in U_h $ \\
\begin{equation}
\label{pbh2}
\tilde{a}_h(u_h,v) = L_h(v),
\end{equation}    

\noindent holds true, where $\forall u \in V+W_h$ and $\forall v \in V+U_h$,\\
\begin{equation}
\label{ahLps2}
\left\{
\begin{array}{l}
\tilde{a}_h(u,v):= \; \displaystyle \sum_{T \in {\cal T}_h} 
\left[ (\nabla \cdot {\cal K} \nabla u, v)_T    
+(\nabla u,{\cal K} \nabla v)_{T} + (u,\nabla \cdot {\cal K} \nabla v)_T \right]\\
L_h(v):=\;-(f,\Pi_h[v]).
\end{array}
\right.
\end{equation}

At a first glance (\ref{ahLps2}) seems to indicate that the velocity ${\bf w}$ does not appear in formulation (\ref{pbh2}). 
Nonetheless ${\bf w}$ remains implicit therein through the definition of space $W_h$.\\
\indent The fact that problem (\ref{pbh2}) has a unique solution can be established quite similarly to problem (\ref{pbh}). 
The convergence results that hold for this method can be proved very much like in the case of the method defined in Section 2. 
The main difference is that it is necessary to require a little more regularity of $\nabla \cdot {\bf w}$, namely,  
that this function lies in $W^{1,\infty}(\Omega)$. 
Apart from this assumption, the results are qualitatively equivalent, in the sense that a priori error estimates completely 
analogous to those of \textbf{Theorem 3.2} and \textbf{Theorem 3.3} apply to problem (\ref{pbh2}) as well. As far as this work 
is concerned, resulting properties among others we have 
not formally established here, are illustrated by means of a numerical example given in the following section.

\section{Numerical experiments}

\hspace{4mm} In this section we present some numerical results obtained with the methods described in Sections 2 and 4 for two 
test problems,  which particularly highlight their behavior. The following nomenclature is used for the different numerical methods having been experimented:
\begin{itemize}
\item Method $A$ - Douglas \& Roberts version in non divergence form of mixed method $RT_0$; 
\item Method $hA$ - Hermite analog of Method $A$ (cf. Section 2);
\item Method $B$ - Douglas \& Roberts version in divergence form of mixed method $RT_0$;
\item Method $hB$ - Hermite analog of Method $B$ (cf. Section 4).
\end{itemize}

\noindent \underline{Test-problem 1:} In these experiments $\Omega$ is the unit square and a manufactured solution $u$ 
is given by $u(x_1,x_2)=(x_1-x_1^2)(x_2-x_2^2)/4$. This together with the choice ${\mathcal K}={\mathcal I}$ and 
${\bf w}=$P\'e$[x_1^2,x_2^2]^{{\mbox t}}/\sqrt{2}$ where P\'e is the P\'eclet number, produces a 
right hand side datum $f$. A sequence of uniform meshes was employed with $2L^2$ triangles, 
for $L=8,16,32,64$, constructed by first subdividing $\Omega$ into $L^2$ equal squares and then each one of these squares into 
two triangles by means of their diagonals parallel to the line $x_1=x_2$. In Figure 1 we display the absolute errors 
in four different respects for increasing values of $L$, of the approximate solutions obtained with methods $A$, $hA$, $B$ and $hB$  
for P\'e$=1$. The notations are self-explanatory, except for $Max |e_u|$ which refers to the maximum of the absolute errors at the centroids of the mesh triangles. Henceforth we call $Max |e_u|$ the \textit{pseudo maximum semi-norm}. More precisely the absolute errors of $u$, $\nabla u$ and $\Delta u = \nabla \cdot {\mathcal K} \nabla u$ measured in the norm of $L^2(\Omega)$ together  
with $Max |e_u|$, are shown in the sub-figures where indicated. In Figures 2 the same kind of results are displayed for P\'e$=100$. \\
From Figures 1 and 2 one can infer that:
\begin{itemize}
  \item Methods $A$ and $hA$ are fairly equivalent to Methods $B$ and $hB$ in all respects for a low P\'eclet number. 
		\item Methods $A$ and $hA$ are superior to Methods $B$ and $hB$ in all respects when the P\'eclet number is not low.
  \item The theoretical results of Section 3 for Method $hA$ were confirmed in the case of both a low and a moderate P\'eclet number.
	\item As the P\'eclet number increases the convergence rate of Method $hB$ in $L^2(\Omega)$ seems to be decreasing from ca. two 
	for P\'e$=1$.  
  \item The numerical convergence rate in the pseudo-maximum semi-norm $Max |e_u|$ is approximately two for 
  all the four methods.
	\item For P\'e$=100$ the mixed methods are a little more accurate than their Hermite counterparts, 
	as far as errors at the triangle centroids are concerned.
	\item For P\'e$=100$ both $\nabla u$ and $\Delta u$ tend to be equally approximated by a mixed method and its 
	Hermite counterpart. 
\end{itemize}
  \begin{figure}[h!]
   \center{\subfigure{\begin{tikzpicture}[scale=0.8]
    \begin{loglogaxis}[title={\textbf{$L^2$ errors of $u$ for P\'e $=1$}},
        width=7cm, height=8cm,
        grid = both,
	tick align=outside,
	tickpos=left,
  	xticklabels={
	  $1/64$,
	  $1/32$,,
	  $1/16$,,,,
	  $1/8$},
  	xtick={0,0.015625,...,1}]
      \addplot[black,mark=*,style=solid] table[x=h,y=h] {dedomena.dat};
      \addplot[black,mark=*,style=densely dotted] table[x=h,y=h2] {dedomena.dat};
      \addplot[mark options={solid,thin,black},gray,mark=square,style=solid,thick] table[x=h,y=SquaUMETa] {dedomena.dat};
      \addplot[mark options={solid,scale=1.3,thick},black,mark=+,style=dashed] table[x=h,y=SquaUMETb] {dedomena.dat};
      \addplot[mark options={solid,scale=1.5,thin,black},gray,mark=x,style=solid,thick] table[x=h,y=SquaUMETha] {dedomena.dat};
      \addplot[mark options={solid,scale=1.0},black,mark=o,style=dashed] table[x=h,y=SquaUMEThb] {dedomena.dat};
    \end{loglogaxis}
\draw[color=black] (5.5,0.0) node [anchor=north west] {$h$};
\draw[color=black] (0.0,6.3) node [anchor=south east] {$\left|\left|e_u\right|\right|_{0,\Omega}$};
\end{tikzpicture}}
   \subfigure{\begin{tikzpicture}[scale=0.8]
    \begin{loglogaxis}[title={\textbf{Maximum errors of $u$ for P\'e $=1$}},
        width=7cm, height=8cm,
        grid = both,
	tick align=outside,
	tickpos=left,
  	xticklabels={
	  $1/64$,
	  $1/32$,,
	  $1/16$,,,,
	  $1/8$},
  	xtick={0,0.015625,...,1}]
      \addplot[mark options={solid},black,mark=*,style=solid] table[x=h,y=h] {dedomena.dat};
      \addplot[mark options={solid},black,mark=*,style=densely dotted] table[x=h,y=h2] {dedomena.dat};
      \addplot[mark options={solid,thin,black},gray,mark=square,style=solid,thick] table[x=h,y=SquaMaxEuMETa] {dedomena.dat};
      \addplot[mark options={solid,scale=1.5,thin,black},gray,mark=x,style=solid,thick] table[x=h,y=SquaMaxEuMETha] {dedomena.dat};
      \addplot[mark options={solid,scale=1.3,thick},black,mark=+,style=dashed] table[x=h,y=SquaMaxEuMETb] {dedomena.dat};
      \addplot[mark options={solid,scale=1.0},black,mark=o,style=dashed] table[x=h,y=SquaMaxEuMEThb] {dedomena.dat};
    \end{loglogaxis}
\draw[color=black] (5.5,0.0) node [anchor=north west] {$h$};
\draw[color=black] (0.0,6.3) node [anchor=south east] {$Max\left|e_{u}\right|$};
\end{tikzpicture}}
   \subfigure{\begin{tikzpicture}[scale=0.8]
    \begin{loglogaxis}[title={\textbf{$L^2$ errors of $\nabla u$ for P\'e $=1$}},
        width=7cm, height=8cm,
        grid = both,
	tick align=outside,
	tickpos=left,
  	xticklabels={
	  $1/64$,
	  $1/32$,,
	  $1/16$,,,,
	  $1/8$},
  	xtick={0,0.015625,...,1}]
      \addplot[mark options={solid},black,mark=*,style=solid] table[x=h,y=h] {dedomena.dat};
      \addplot[mark options={solid},black,mark=*,style=densely dotted] table[x=h,y=h2] {dedomena.dat};
      \addplot[mark options={solid,thin,black,scale=1.3},gray,mark=square,style=solid,thick] table[x=h,y=SquaGraduMETa] {dedomena.dat};
      \addplot[mark options={solid,scale=1.5,thin,black},gray,mark=x,style=solid,thick] table[x=h,y=SquaGraduMETha] {dedomena.dat};
      \addplot[mark options={solid,scale=1.3,thick},black,mark=+,style=dashed] table[x=h,y=SquaGraduMETb] {dedomena.dat};
      \addplot[mark options={solid,scale=1.0},black,mark=o,style=dashed] table[x=h,y=SquaGraduMEThb] {dedomena.dat};
    \end{loglogaxis}
\draw[color=black] (5.5,0.0) node [anchor=north west] {$h$};
\draw[color=black] (0.0,6.3) node [anchor=south east] {$\left|\left|e_{\nabla u}\right|\right|_{0,\Omega}$};
\end{tikzpicture}}
   \subfigure{\begin{tikzpicture}[scale=0.8]
    \begin{loglogaxis}[title={\textbf{$L^2$ errors of $\Delta u$ for P\'e $=1$}},
        width=7cm, height=8cm,
        grid = both,
	tick align=outside,
	tickpos=left,
  	xticklabels={
	  $1/64$,
	  $1/32$,,
	  $1/16$,,,,
	  $1/8$},
  	xtick={0,0.015625,...,1}]
      \addplot[mark options={solid},black,mark=*,style=solid] table[x=h,y=h] {dedomena.dat};
      \addplot[mark options={solid},black,mark=*,style=densely dotted] table[x=h,y=h2] {dedomena.dat};
      \addplot[mark options={solid,thin,black,scale=1.3},gray,mark=square,style=solid,thick] table[x=h,y=SquaDeltuMETa] {dedomena.dat};
      \addplot[mark options={solid,scale=1.3,thick},black,mark=+,style=dashed] table[x=h,y=SquaDeltuMETb] {dedomena.dat};
      \addplot[mark options={solid,scale=1.5,thin,black},gray,mark=x,style=solid,thick] table[x=h,y=SquaDeltuMETha] {dedomena.dat};
      \addplot[mark options={solid,scale=1.0},black,mark=o,style=dashed] table[x=h,y=SquaDeltuMEThb] {dedomena.dat};
    \end{loglogaxis}
\draw[color=black] (5.5,0.0) node [anchor=north west] {$h$};
\draw[color=black] (0.0,6.3) node [anchor=south east] {$\left|\left|e_{\Delta u}\right|\right|_{0,\Omega}$};
\end{tikzpicture}}
   \subfigure{\begin{tikzpicture}
  \coordinate (ih) at (0,.5);
  \coordinate [label=right:$h$] (h) at (.6,.5);
  \coordinate (markh) at (.3,.5);
  \coordinate (ih2) at (1.5,.5);
  \coordinate [label=right:$h^2$] (h2) at (2.1,.5);
  \coordinate (markh2) at (1.8,.5);
  \coordinate (iA) at (3.0,.5);
  \coordinate [label=right:Method $A$] (A) at (3.6,.5);
  \coordinate (markA) at (3.3,.5);
  \coordinate (iB) at (5.8,.5);
  \coordinate [label=right:Method $B$] (B) at (6.4,.5);
  \coordinate (markB) at (6.1,.5);
  \coordinate (ihA) at (8.6,.5);
  \coordinate [label=right:Method $hA$] (hA) at (9.2,.5);
  \coordinate (markhA) at (8.9,.5);
  \coordinate (ihB) at (11.6,.5);
  \coordinate [label=right:Method $hB$] (hB) at (12.2,.5);
  \coordinate (markhB) at (11.9,.5);
  \draw (-.3,0) rectangle (14.3,1);
  \draw (ih) -- plot[mark=*] coordinates {(markh)} -- (h);
  \draw[densely dotted] (ih2) -- plot[mark=*] coordinates {(markh2)} -- (h2);
  \draw[gray,thick] (iA) -- plot[mark=square,mark options={solid,thin,black}] coordinates {(markA)} -- (A);
  \draw[dashed] (iB) -- plot[mark=+,mark options={solid,scale=1.3,thick}] coordinates {(markB)} -- (B);
  \draw[gray,thick] (ihA) -- plot[mark=x,mark options={solid,scale=1.5,thin,black}] coordinates {(markhA)} -- (hA);
  \draw[dashed] (ihB) -- plot[mark=o,mark options={solid,scale=1.0}] coordinates {(markhB)} -- (hB);
\end{tikzpicture}}}
  \caption{Results for Test-problem 1 with P\'e $=1$}
  \label{fig:SquaPe1}
  \end{figure}
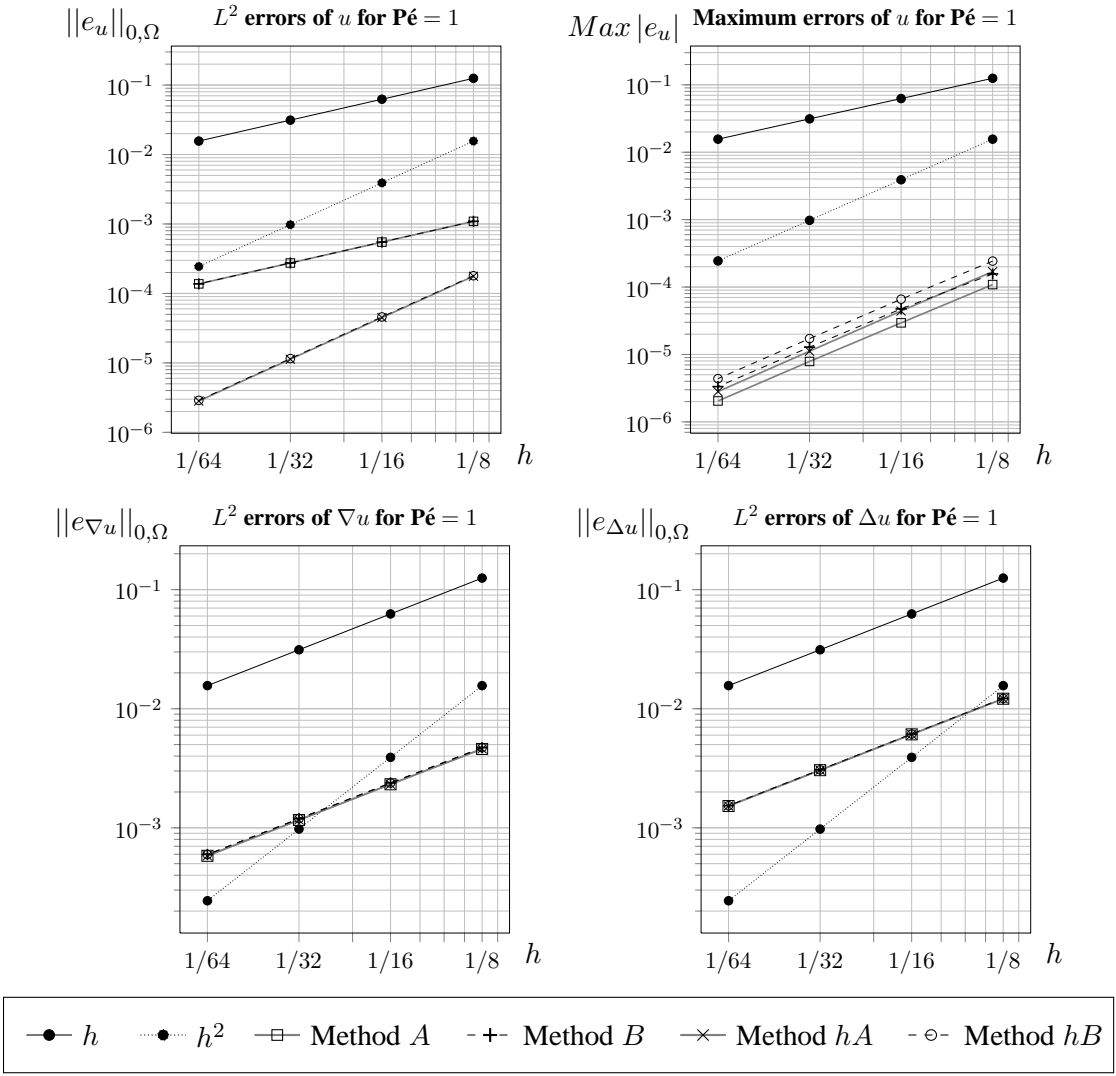
  \begin{figure}[h!]
   \center{\subfigure{\begin{tikzpicture}[scale=.8]
    \begin{loglogaxis}[title={\textbf{$L^2$ errors of $u$ for P\'e $=100$}},
        width=7cm, height=8cm,
        grid = both,
	tick align=outside,
	tickpos=left,
  	xticklabels={
	  $1/64$,
	  $1/32$,,
	  $1/16$,,,,
	  $1/8$},
  	xtick={0,0.015625,...,1}]
      \addplot[black,mark=*,style=solid] table[x=h,y=h] {dedomena.dat};
      \addplot[black,mark=*,style=densely dotted] table[x=h,y=h2] {dedomena.dat};
      \addplot[black,mark=square,style=solid] table[x=h,y=S100UMETa] {dedomena.dat};
      \addplot[mark options={solid,scale=1.3,thick},black,mark=+,style=dashed] table[x=h,y=S100UMETb] {dedomena.dat};
      \addplot[mark options={scale=1.5},black,mark=x,style=solid] table[x=h,y=S100UMETha] {dedomena.dat};
      \addplot[mark options={solid,scale=1.0},black,mark=o,style=dashed] table[x=h,y=S100UMEThb] {dedomena.dat};
    \end{loglogaxis}
\draw[color=black] (5.5,0.0) node [anchor=north west] {$h$};
\draw[color=black] (0.0,6.3) node [anchor=south east] {$\left|\left|e_u\right|\right|_{0,\Omega}$};
\end{tikzpicture}}
   \subfigure{\begin{tikzpicture}[scale=.8]
    \begin{loglogaxis}[title={\textbf{Maximum errors of $u$ for P\'e $=100$}},
        width=7cm, height=8cm,
        grid = both,
	tick align=outside,
	tickpos=left,
  	xticklabels={
	  $1/64$,
	  $1/32$,,
	  $1/16$,,,,
	  $1/8$},
  	xtick={0,0.015625,...,1}]
      \addplot[mark options={solid},black,mark=*,style=solid] table[x=h,y=h] {dedomena.dat};
      \addplot[mark options={solid},black,mark=*,style=densely dotted] table[x=h,y=h2] {dedomena.dat};
      \addplot[black,mark=square,style=solid] table[x=h,y=S100MaxEuMETa] {dedomena.dat};
      \addplot[mark options={solid,scale=1.3,thick},black,mark=+,style=dashed] table[x=h,y=S100MaxEuMETb] {dedomena.dat};
      \addplot[mark options={scale=1.5},black,mark=x,style=solid] table[x=h,y=S100MaxEuMETha] {dedomena.dat};
      \addplot[mark options={solid,scale=1.0},black,mark=o,style=dashed] table[x=h,y=S100MaxEuMEThb] {dedomena.dat};
    \end{loglogaxis}
\draw[color=black] (5.5,0.0) node [anchor=north west] {$h$};
\draw[color=black] (0.0,6.3) node [anchor=south east] {$Max\left|e_{u}\right|$};
\end{tikzpicture}}
   \subfigure{\begin{tikzpicture}[scale=.8]
    \begin{loglogaxis}[title={\textbf{$L^2$ errors of $\nabla u$ for P\'e $=100$}},
        width=7cm, height=8cm,
        grid = both,
	tick align=outside,
	tickpos=left,
  	xticklabels={
	  $1/64$,
	  $1/32$,,
	  $1/16$,,,,
	  $1/8$},
  	xtick={0,0.015625,...,1}]
      \addplot[mark options={solid,black,thin},gray,mark=*,style=solid,thick] table[x=h,y=h] {dedomena.dat};
      \addplot[mark options={solid},black,mark=*,style=densely dotted] table[x=h,y=h2] {dedomena.dat};
      \addplot[black,mark=square,style=solid] table[x=h,y=S100GraduMETa] {dedomena.dat};
      \addplot[mark options={solid,scale=1.3,thick},black,mark=+,style=dashed] table[x=h,y=S100GraduMETb] {dedomena.dat};
      \addplot[mark options={scale=1.5},black,mark=x,style=solid] table[x=h,y=S100GraduMETha] {dedomena.dat};
      \addplot[mark options={solid,scale=1.0},black,mark=o,style=dashed] table[x=h,y=S100GraduMEThb] {dedomena.dat};
    \end{loglogaxis}
\draw[color=black] (5.5,0.0) node [anchor=north west] {$h$};
\draw[color=black] (0.0,6.3) node [anchor=south east] {$\left|\left|e_{\nabla u}\right|\right|_{0,\Omega}$};
\end{tikzpicture}}
   \subfigure{\begin{tikzpicture}[scale=.8]
    \begin{loglogaxis}[title={\textbf{$L^2$ errors of $\Delta u$ for P\'e $=100$}},
        width=7cm, height=8cm,
        grid = both,
	tick align=outside,
	tickpos=left,
  	xticklabels={
	 $1/64$,
	  $1/32$,,
	  $1/16$,,,,
	  $1/8$},
  	xtick={0,0.015625,...,1}]
      \addplot[mark options={solid},black,mark=*,style=solid] table[x=h,y=h] {dedomena.dat};
      \addplot[mark options={solid},black,mark=*,style=densely dotted] table[x=h,y=h2] {dedomena.dat};
      \addplot[black,mark=square,style=solid] table[x=h,y=S100DeltuMETa] {dedomena.dat};
      \addplot[mark options={solid,scale=1.3,thick},black,mark=+,style=dashed] table[x=h,y=S100DeltuMETb] {dedomena.dat};
      \addplot[mark options={scale=1.5,thick,gray},black,mark=x,style=solid] table[x=h,y=S100DeltuMETha] {dedomena.dat};
      \addplot[mark options={solid,scale=1.0},black,mark=o,style=dashed] table[x=h,y=S100DeltuMEThb] {dedomena.dat};
    \end{loglogaxis}
\draw[color=black] (5.5,0.0) node [anchor=north west] {$h$};
\draw[color=black] (0.0,6.3) node [anchor=south east] {$\left|\left|e_{\Delta u}\right|\right|_{0,\Omega}$};
\end{tikzpicture}}
   \subfigure{\begin{tikzpicture}
  \coordinate (ih) at (0,.5);
  \coordinate [label=right:$h$] (h) at (.6,.5);
  \coordinate (markh) at (.3,.5);
  \coordinate (ih2) at (1.5,.5);
  \coordinate [label=right:$h^2$] (h2) at (2.1,.5);
  \coordinate (markh2) at (1.8,.5);
  \coordinate (iA) at (3.0,.5);
  \coordinate [label=right:Method $A$] (A) at (3.6,.5);
  \coordinate (markA) at (3.3,.5);
  \coordinate (iB) at (5.8,.5);
  \coordinate [label=right:Method $B$] (B) at (6.4,.5);
  \coordinate (markB) at (6.1,.5);
  \coordinate (ihA) at (8.6,.5);
  \coordinate [label=right:Method $hA$] (hA) at (9.2,.5);
  \coordinate (markhA) at (8.9,.5);
  \coordinate (ihB) at (11.6,.5);
  \coordinate [label=right:Method $hB$] (hB) at (12.2,.5);
  \coordinate (markhB) at (11.9,.5);
  \draw (-.3,0) rectangle (14.3,1);
  \draw (ih) -- plot[mark=*] coordinates {(markh)} -- (h);
  \draw[densely dotted] (ih2) -- plot[mark=*] coordinates {(markh2)} -- (h2);
  \draw[gray,thick] (iA) -- plot[mark=square,mark options={solid,thin,black}] coordinates {(markA)} -- (A);
  \draw[dashed] (iB) -- plot[mark=+,mark options={solid,scale=1.3,thick}] coordinates {(markB)} -- (B);
  \draw[gray,thick] (ihA) -- plot[mark=x,mark options={solid,scale=1.5,thin,black}] coordinates {(markhA)} -- (hA);
  \draw[dashed] (ihB) -- plot[mark=o,mark options={solid,scale=1.0}] coordinates {(markhB)} -- (hB);
\end{tikzpicture}}}
  \caption{Results for Test-problem 1 with P\'e $=100$}
  \label{fig:SquaPe100}
  \end{figure}
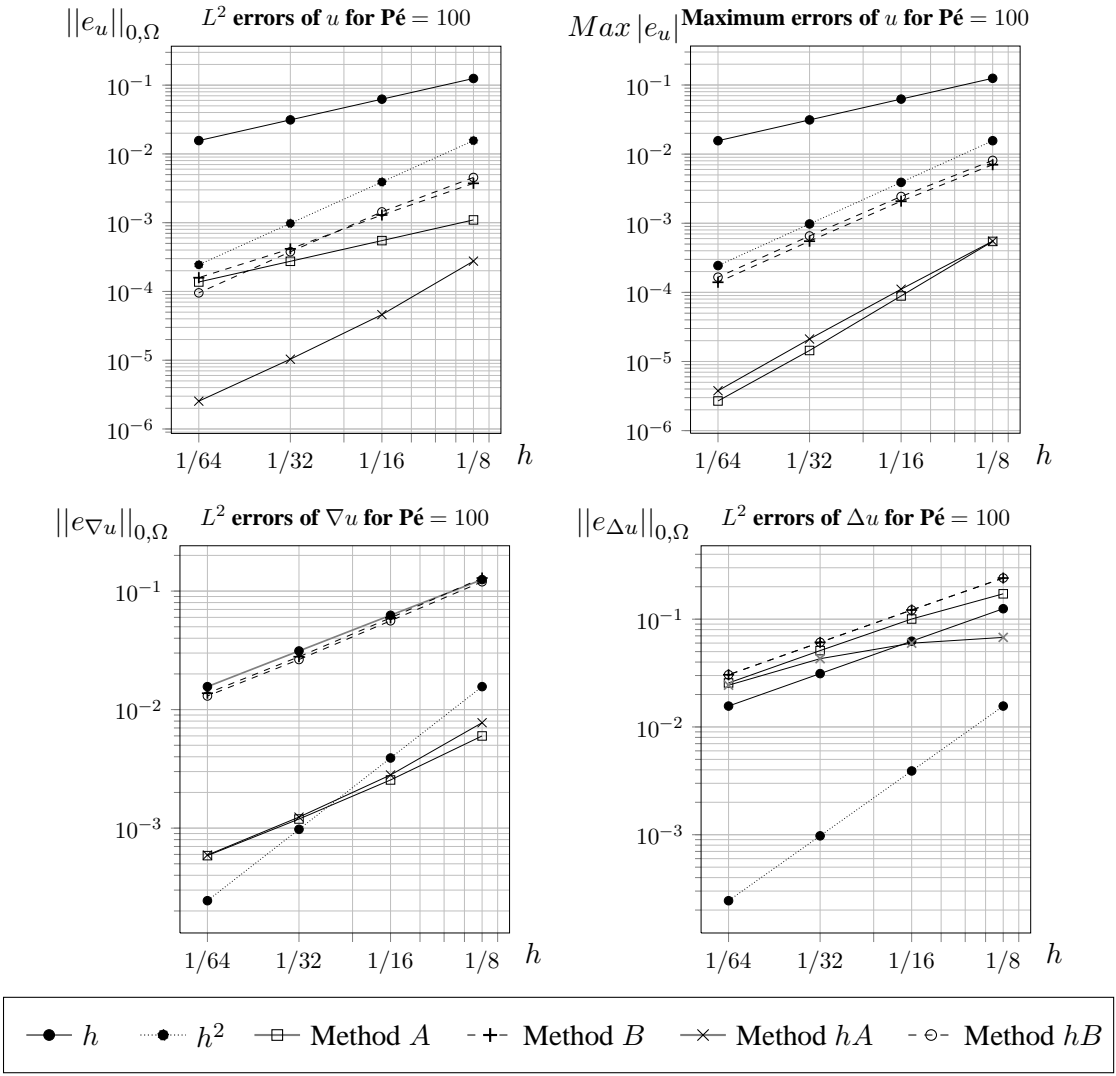
\clearpage
The numerical results for Methods $B$ and $hB$ deteriorated substantially as we switched to higher P\'eclet numbers, which was partially the case of Method $hA$, while most of the results obtained with Method $A$ remained quite reasonable. 
Taking $L=64$ we illustrate the behavior of Methods $A$ and $hA$ in Tables 1 and 2, respectively, for increasing P\'eclet numbers. 
More precisely we took P\'e$=10^{2k}$ for $k=0,1,2,3$, for which we display the absolute errors of $u$, $\nabla u$, $\Delta u$ measured in the norm of $L^2(\Omega)$ and the absolute error of $u$ measured in the pseudo maximum semi-norm $Max|e_u|$.\\
\begin{table}[h!]
\begin{center}
\begin{tabular}{|c|c|c|c|c|}
\hline P\'e & $\left|\left| e_u \right|\right|_{0,\Omega}$ & $\left|\left| e_{\nabla u} \right|\right|_{0,\Omega}$  
& $\left|\left| e_{\Delta u} \right|\right|_{0,\Omega}$ & $Max \left| e_{u} \right|$ \\
\hline \hline
\footnotesize{1}
 & \footnotesize{0.13723841E-03}
 & \footnotesize{0.58218263E-03}
 & \footnotesize{0.15249263E-02}
 & \footnotesize{0.20428256E-05}
\\
\hline
\footnotesize{100}
 & \footnotesize{0.13724039E-03}
 & \footnotesize{0.58595099E-03}
 & \footnotesize{0.25587661E-01}
 & \footnotesize{0.26841141E-05}
\\
\hline
\footnotesize{10000}
 & \footnotesize{0.18370239E-02}
 & \footnotesize{0.51526514E+00}
 & \footnotesize{0.15093095E+03}
 & \footnotesize{0.52745011E-01}
\\
\hline
\footnotesize{1000000}
 & \footnotesize{0.20738981E-03}
 & \footnotesize{0.57979017E-01}
 & \footnotesize{0.22150639E+02}
 & \footnotesize{0.12137294E-02}
\\
\hline
\end{tabular} \label{tabelaSQUAa64}
\caption{Absolute errors for Test-problem 1 solved by Method $A$}
\end{center}
\end{table}
\begin{table}[h!]
\begin{center}
\begin{tabular}{|c|c|c|c|c|}
\hline P\'e & $\left|\left| e_u \right|\right|_{0,\Omega}$ & $\left|\left| e_{\nabla u} \right|\right|_{0,\Omega}$ 
& $\left|\left| e_{\Delta u} \right|\right|_{0,\Omega}$ & $Max \left| e_{u} \right|$ \\
\hline \hline
\footnotesize{1}
 & \footnotesize{0.28250216E-05}
 & \footnotesize{0.58219418E-03}
 & \footnotesize{0.15249297E-02}
 & \footnotesize{0.28130033E-05}
\\
\hline
\footnotesize{100}
 & \footnotesize{0.25386722E-05}
 & \footnotesize{0.59256341E-03}
 & \footnotesize{0.24402972E-01}
 & \footnotesize{0.37752993E-05}
\\
\hline
\footnotesize{10000}
 & \footnotesize{0.10089341E+04}
 & \footnotesize{0.26546637E+06}
 & \footnotesize{0.75214576E+07}
 & \footnotesize{0.57326794E+04}
\\
\hline
\footnotesize{1000000}
 & \footnotesize{0.13681695E+00}
 & \footnotesize{0.26569231E+02}
 & \footnotesize{0.66544379E+02}
 & \footnotesize{0.27070200E+01}
\\
\hline
\end{tabular} \label{tabelaSQUAha64}
\caption{Absolute errors for Test-problem 1 solved by Method $hA$}
\end{center}
\end{table}
\noindent From Tables 1 and 2 we conclude that (mixed) Method A is more stable than its Hermite counterpart $hA$, as the P\'eclet  number increases. \\

\noindent \underline{Test-problem 2:}  In order to observe the behavior of the four methods being checked, in the presence of 
a curved boundary, in this test-problem the domain is a disk with unit radius. The manufactured solution $u$ 
is given by $u(x_1,x_2)=(1-x_1^2-x_2^2)/4$. Taking again ${\mathcal K}={\mathcal I}$, 
the right hand side function $f=1-(x_1^2+x_2^2)/2$ corresponds to a convective velocity ${\bf w}=$P\'e$[x_1,x_2]^{{\mbox t}}$. 
For symmetry reasons the computational domain $\Omega$ is only the quarter of disk given by $x_1>0$ and $x_2>0$. 
A sequence of quasi-uniform meshes with $2L^2$ triangles was employed 
for $L=8,16,32,64$, constructed by mapping the meshes of Test-problem 1 into the actual meshes of $\Omega$ using 
the transformation of cartesian into polar coordinates in the way described in \cite{CMA}. We denote by $\Omega_h$ the approximation 
of $\Omega$ consisting of the union of the triangles in ${\mathcal T}_h$. In Figure 3 we display the absolute errors 
in four different respects for increasing values of $L$, of the approximate solutions obtained with methods $A$, $hA$, $B$ and $hB$  
for P\'e$=1$. The displayed errors and corresponding notations are the same as in Figures 1 and 2.
More precisely the absolute errors of $u$, $\nabla u$ and $\Delta u = \nabla \cdot {\mathcal K} \nabla u$ 
measured in the norm of $L^2(\Omega_h)$ together  with $Max |e_u|$, are shown in four sub-figures. \\
From Figure 3 we infer that:
\begin{itemize}
  \item Methods $A$ and $hA$ are superior to Methods $B$ and $hB$ in all respects, except in the approximation of (constant) 
	$\Delta u$, which is almost exactly approximated by all the four methods.
		\item Methods $A$ and $hA$ do not seem to be affected by the curved boundary approximation by polygons, while this seems to be case of Methods $B$ and $hB$.
  \item Method $A$ and $hA$ approximate both $\nabla u$ and $\Delta u$ to machine precision; this is an expected behavior  
	since both functions in this test-problem can be exactly represented by the same underlying incomplete linear and 
	constant interpolation for both methods. 
	\item The approximations of $u$ by Method $hA$ converge as an $O(h^2)$ in $L^2(\Omega_h)$, while those computed by Method $A$ 
	converge as an $O(h)$, i.e. the best we can hope for.
  \item The numerical convergence rate in the pseudo-maximum semi-norm $Max |e_u|$ is approximately two for 
 both Method $A$ and Method $hA$ with an advantage of the former over the latter in terms of accuracy. \\
\end{itemize}

  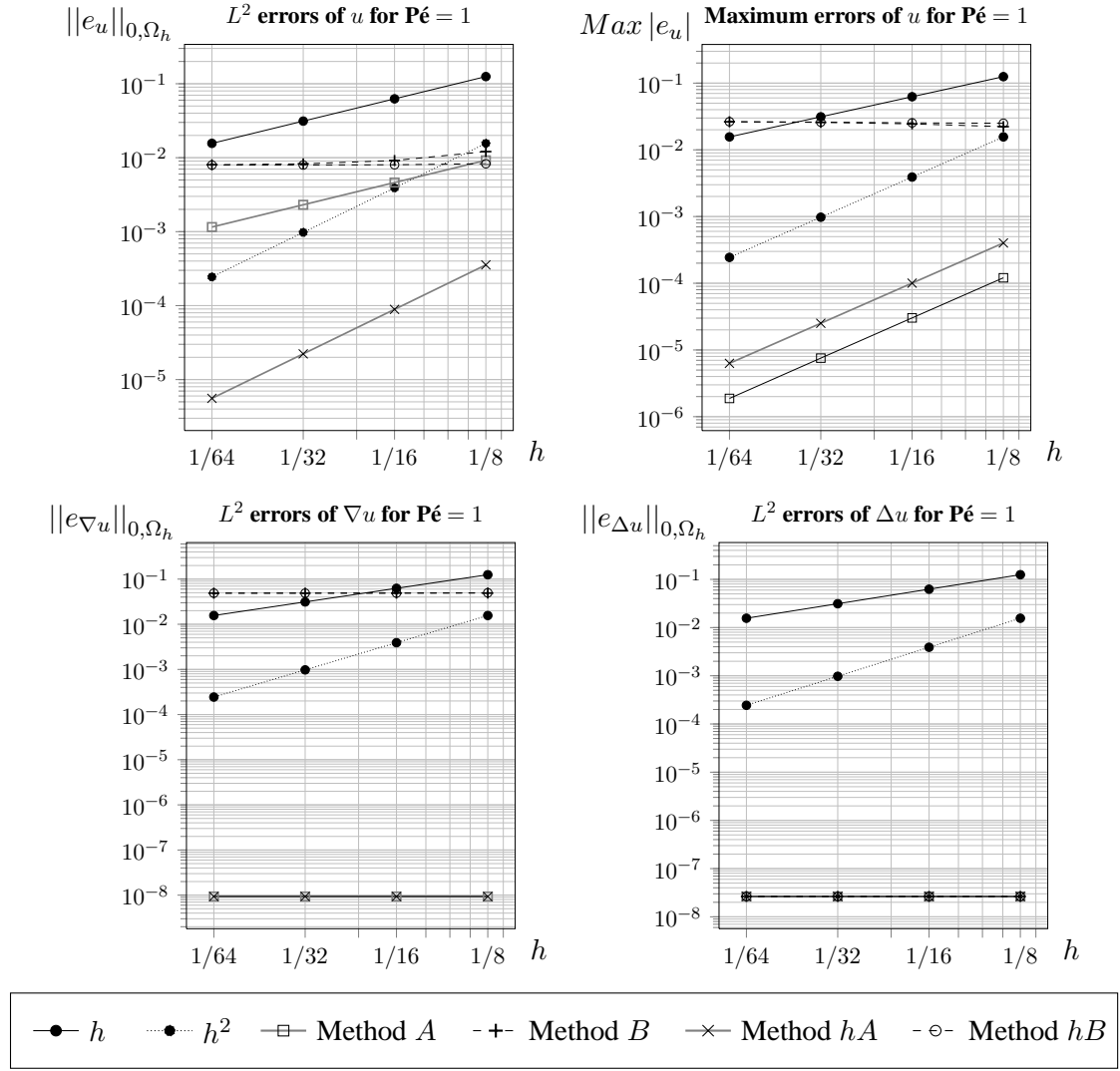
\begin{figure}[h!]
   \center{\subfigure{\begin{tikzpicture}[scale=.8]
    \begin{loglogaxis}[title={\textbf{$L^2$ errors of $u$ for P\'e $=1$}},
        width=7cm, height=8cm,
        grid = both,
	tick align=outside,
	tickpos=left,
  	xticklabels={
	  $1/64$,
	  $1/32$,,
	  $1/16$,,,,
	  $1/8$},
  	xtick={0,0.015625,...,1}]
      \addplot[black,mark=*,style=solid] table[x=h,y=h] {dedomena.dat};
      \addplot[black,mark=*,style=densely dotted] table[x=h,y=h2] {dedomena.dat};
      \addplot[mark options={solid},gray,mark=square,style=solid,thick] table[x=h,y=DiskUMETa] {dedomena.dat};
      \addplot[mark options={solid,scale=1.3,thick},black,mark=+,style=dashed] table[x=h,y=DiskUMETb] {dedomena.dat};
      \addplot[mark options={solid,scale=1.5,thin,black},gray,mark=x,style=solid,thick] table[x=h,y=DiskUMETha] {dedomena.dat};
      \addplot[mark options={solid,scale=1.0},black,mark=o,style=dashed] table[x=h,y=DiskUMEThb] {dedomena.dat};
    \end{loglogaxis}
\draw[color=black] (5.5,0.0) node [anchor=north west] {$h$};
\draw[color=black] (0.0,6.3) node [anchor=south east] {$\left|\left|e_u\right|\right|_{0,\Omega_h}$};
\end{tikzpicture}}
   \subfigure{\begin{tikzpicture}[scale=.8]
    \begin{loglogaxis}[title={\textbf{Maximum errors of $u$ for P\'e $=1$}},
        width=7cm, height=8cm,
        grid = both,
	tick align=outside,
	tickpos=left,
  	xticklabels={
	  $1/64$,
	  $1/32$,,
	  $1/16$,,,,
	  $1/8$},
  	xtick={0,0.015625,...,1}]
      \addplot[mark options={solid},black,mark=*,style=solid] table[x=h,y=h] {dedomena.dat};
      \addplot[mark options={solid},black,mark=*,style=densely dotted] table[x=h,y=h2] {dedomena.dat};
      \addplot[mark options={solid},black,mark=square,style=solid] table[x=h,y=DiskMaxEuMETa] {dedomena.dat};
      \addplot[mark options={solid,scale=1.3,thick},black,mark=+,style=dashed] table[x=h,y=DiskMaxEuMETb] {dedomena.dat};
      \addplot[mark options={solid,scale=1.5,thin,black},gray,mark=x,style=solid,thick] table[x=h,y=DiskMaxEuMETha] {dedomena.dat};
      \addplot[mark options={solid,scale=1.0},black,mark=o,style=dashed] table[x=h,y=DiskMaxEuMEThb] {dedomena.dat};
    \end{loglogaxis}
\draw[color=black] (5.5,0.0) node [anchor=north west] {$h$};
\draw[color=black] (0.0,6.3) node [anchor=south east] {$Max\left|e_{u}\right|$};
\end{tikzpicture}}
   \subfigure{\begin{tikzpicture}[scale=.8]
    \begin{loglogaxis}[title={\textbf{$L^2$ errors of $\nabla u$ for P\'e $=1$}},
        width=7cm, height=8cm,
        grid = both,
	tick align=outside,
	tickpos=left,
  	xticklabels={
	  $1/64$,
	  $1/32$,,
	  $1/16$,,,,
	  $1/8$},
  	xtick={0,0.015625,...,1}]
      \addplot[mark options={solid},black,mark=*,style=solid] table[x=h,y=h] {dedomena.dat};
      \addplot[mark options={solid},black,mark=*,style=densely dotted] table[x=h,y=h2] {dedomena.dat};
      \addplot[mark options={solid},gray,thick,mark=square,style=solid] table[x=h,y=DiskGraduMETa] {dedomena.dat};
      \addplot[mark options={solid,scale=1.3,thick},black,mark=+,style=dashed] table[x=h,y=DiskGraduMETb] {dedomena.dat};
      \addplot[mark options={solid,scale=1.5,thin,black},black,mark=x,style=solid] table[x=h,y=DiskGraduMETha] {dedomena.dat};
      \addplot[mark options={solid,scale=1.0},black,mark=o,style=dashed] table[x=h,y=DiskGraduMEThb] {dedomena.dat};
    \end{loglogaxis}
\draw[color=black] (5.5,0.0) node [anchor=north west] {$h$};
\draw[color=black] (0.0,6.3) node [anchor=south east] {$\left|\left|e_{\nabla u}\right|\right|_{0,\Omega_h}$};
\end{tikzpicture}}
   \subfigure{\begin{tikzpicture}[scale=.8]
    \begin{loglogaxis}[title={\textbf{$L^2$ errors of $\Delta u$ for P\'e $=1$}},
        width=7cm, height=8cm,
        grid = both,
	tick align=outside,
	tickpos=left,
  	xticklabels={
	  $1/64$,
	  $1/32$,,
	  $1/16$,,,,
	  $1/8$},
  	xtick={0,0.015625,...,1}]
      \addplot[mark options={solid},black,mark=*,style=solid] table[x=h,y=h] {dedomena.dat};
      \addplot[mark options={solid},black,mark=*,style=densely dotted] table[x=h,y=h2] {dedomena.dat};
      \addplot[mark options={solid},gray,thick,mark=square,style=solid] table[x=h,y=DiskDeltuMETa] {dedomena.dat};
      \addplot[mark options={solid,scale=1.3,thick},black,mark=+,style=dashed] table[x=h,y=DiskDeltuMETb] {dedomena.dat};
      \addplot[mark options={solid,scale=1.5,thin,black},black,mark=x,style=solid] table[x=h,y=DiskDeltuMETha] {dedomena.dat};
      \addplot[mark options={solid,scale=1.0},black,mark=o,style=dashed] table[x=h,y=DiskDeltuMEThb] {dedomena.dat};
    \end{loglogaxis}
\draw[color=black] (5.5,0.0) node [anchor=north west] {$h$};
\draw[color=black] (0.0,6.3) node [anchor=south east] {$\left|\left|e_{\Delta u}\right|\right|_{0,\Omega_h}$};
\end{tikzpicture}}
   \subfigure{\begin{tikzpicture}
  \coordinate (ih) at (0,.5);
  \coordinate [label=right:$h$] (h) at (.6,.5);
  \coordinate (markh) at (.3,.5);
  \coordinate (ih2) at (1.5,.5);
  \coordinate [label=right:$h^2$] (h2) at (2.1,.5);
  \coordinate (markh2) at (1.8,.5);
  \coordinate (iA) at (3.0,.5);
  \coordinate [label=right:Method $A$] (A) at (3.6,.5);
  \coordinate (markA) at (3.3,.5);
  \coordinate (iB) at (5.8,.5);
  \coordinate [label=right:Method $B$] (B) at (6.4,.5);
  \coordinate (markB) at (6.1,.5);
  \coordinate (ihA) at (8.6,.5);
  \coordinate [label=right:Method $hA$] (hA) at (9.2,.5);
  \coordinate (markhA) at (8.9,.5);
  \coordinate (ihB) at (11.6,.5);
  \coordinate [label=right:Method $hB$] (hB) at (12.2,.5);
  \coordinate (markhB) at (11.9,.5);
  \draw (-.3,0) rectangle (14.3,1);
  \draw (ih) -- plot[mark=*] coordinates {(markh)} -- (h);
  \draw[densely dotted] (ih2) -- plot[mark=*] coordinates {(markh2)} -- (h2);
  \draw[gray,thick] (iA) -- plot[mark=square,mark options={solid,thin,black}] coordinates {(markA)} -- (A);
  \draw[dashed] (iB) -- plot[mark=+,mark options={solid,scale=1.3,thick}] coordinates {(markB)} -- (B);
  \draw[gray,thick] (ihA) -- plot[mark=x,mark options={solid,scale=1.5,thin,black}] coordinates {(markhA)} -- (hA);
  \draw[dashed] (ihB) -- plot[mark=o,mark options={solid,scale=1.0}] coordinates {(markhB)} -- (hB);
\end{tikzpicture}}}
  \caption{Results for Test-problem 2 with P\'e $=1$}
  \label{fig:DiskPe1}
  \end{figure}

Akin to Test-problem 1, we checked the behavior of Methods $A$ and $hA$ as the P\'eclet number increases. Here again we took 
$L=64$ and P\'e$=10^{2k}$ for $k=0,1,2,3$. The resulting errors measured in the same manner as in Tables 1 and 2 are displayed 
in Table 3 for Method $A$ and in Table 4 for Method $hA$. \\
\begin{table}[h!]
\begin{center}
\begin{tabular}{|c|c|c|c|c|}
\hline P\'e & $\left|\left| e_u \right|\right|_{0,\Omega_h}$ & $\left|\left| e_{\nabla u} \right|\right|_{0,\Omega_h}$ 
 & $\left|\left| e_{\Delta u} \right|\right|_{0,\Omega_h}$ & $Max \left| e_{u} \right|$ \\
\hline \hline
\footnotesize{1}
 & \footnotesize{0.11539009E-02}
 & \footnotesize{0.93376835E-08}
 & \footnotesize{0.26411289E-07}
 & \footnotesize{0.18777283E-05}
\\
\hline
\footnotesize{100}
 & \footnotesize{0.11539009E-02}
 & \footnotesize{0.93376819E-08}
 & \footnotesize{0.26411285E-07}
 & \footnotesize{0.18777283E-05}
\\
\hline
\footnotesize{10000}
 & \footnotesize{0.11539009E-02}
 & \footnotesize{0.93377017E-08}
 & \footnotesize{0.26411346E-07}
 & \footnotesize{0.18777283E-05}
\\
\hline
\footnotesize{1000000}
 & \footnotesize{0.11539009E-02}
 & \footnotesize{0.96367921E-08}
 & \footnotesize{0.34526072E-07}
 & \footnotesize{0.18774905E-05}
\\
\hline
\end{tabular} \label{tabelaDISKa64}
\caption{Absolute errors for Test-problem 2 solved by Method $A$}
\end{center}
\end{table}
\begin{table}[h!]
\begin{center}
\begin{tabular}{|c|c|c|c|c|}
\hline P\'e & $\left|\left| e_u \right|\right|_{0,\Omega_h}$  
& $\left|\left| e_{\nabla u} \right|\right|_{0,\Omega_h}$ & $\left|\left| e_{\Delta u} \right|\right|_{0,\Omega_h}$ & $Max \left| e_{u} \right|$ \\
\hline \hline
\footnotesize{1}
 & \footnotesize{0.55709298E-05}
 & \footnotesize{0.93376800E-08}
 & \footnotesize{0.26411287E-07}
 & \footnotesize{0.63028673E-05}
\\
\hline
\footnotesize{100}
 & \footnotesize{0.55708924E-05}
 & \footnotesize{0.92319279E-08}
 & \footnotesize{0.26157933E-07}
 & \footnotesize{0.63028012E-05}
\\
\hline
\footnotesize{10000}
 & \footnotesize{0.55694727E-05}
 & \footnotesize{0.52231986E-08}
 & \footnotesize{0.15104352E-07}
 & \footnotesize{0.62996183E-05}
\\
\hline
\footnotesize{1000000}
 & \footnotesize{0.55709889E-05}
 & \footnotesize{0.93376916E-08}
 & \footnotesize{0.26406856E-07}
 & \footnotesize{0.63029895E-05}
\\
\hline
\end{tabular} \label{tabelaDISKha64}
\caption{Absolute errors for Test-problem 2 solved by Method $hA$}
\end{center}
\end{table}
\noindent From Tables 3 and 4 we observe that both Method $A$ and Method $hA$ are accurate to machine precision, irrespective of the P\'eclet number, as far as the approximations of $\nabla u$ and $\Delta u$ are concerned. The approximations of $u$ in $L^2(\Omega_h)$ 
and at the triangle centroids do not seem to be affected by the P\'eclet number either in this test-problem for both methods. 
In the former sense Method $hA$ is much more accurate than Method $A$ as expected, while in the latter sense Method $A$ is slightly  more precise than Method $hA$. Notice that this test-problem is a little peculiar, since the exact solution is a quadratic 
function, whose gradient can be exactly represented by the gradient of the underlying interpolating functions. Actually 
this also happens to the approximation of the function itself by Method $hA$, but in this case other sources of errors came into play, 
such as numerical integration (see also Remark 3 hereafter).  
\\

\section{Concluding remarks}

We conclude this work with a few remarks.

\begin{remark}
By means of similar test-problems at low to moderate P\'eclet numbers, it was shown in \cite{ICNAAM2013} and \cite{HEFAT} 
that the Hermite methods work as well as the corresponding Douglas and Roberts extensions of the $RT_0$ element, as far as the 
fluxes are concerned. However the former behave much better in terms of the error of the primal variable in $L^2(\Omega)$, 
as expected. \rule{2mm}{2mm}
\end{remark} 

\begin{remark}
According to the theoretical results derived in this work, numerical convergence in case the P\'eclet number is high could only be 
observed if meshes much finer than those used in Test-problem 1 and in \cite{ICNAAM2013} and \cite{HEFAT} were used. However 
running tests with such meshes may become unrealistic. Therefore the authors intend to study modifications of the variational 
formulations employed in this work, in order to obtain stable solutions within acceptable accuracy, even in the case of high 
P\'eclet numbers, without resorting to excessive mesh refinement. The work of Park and Kim \cite{Kim} for $RT_0$ discretizations 
could be a inspiring one in this connection. 
\rule{2mm}{2mm} 
\end{remark}

\begin{remark}
The bilinear forms and the linear form $L_h$ considered in this work do not really reduce to those in \cite{CAM2013} in the case 
where $\bw \equiv \bf 0$. This is because somehow we wanted to incorporate numerical quadrature to the variational formulations in 
use, which is mandatory if $f$ is not easy to integrate. However in this case we should rather take $L_h(v):=(f_h,v)$ for a 
suitable $f_h$ defined through point values of $f$ only. Assuming for instance that $f \in H^2(\Omega)$, $f_h$ can be chosen to be 
a piecewise linear interpolate of $f$ in every $T \in {\mathcal T}_h$. Second order convergence results in $L^2(\Omega)$ can still be proven to 
hold for such a choice, using the well-known analysis of variational crimes \cite{Strang}. The same qualitative results can also be 
obtained by using a suitable quadrature formula to compute the integral of the function $g := f v$ in every element of the mesh. 
For more details we also refer to \cite{Strang}, or to many other text books on the finite element method.
\rule{2mm}{2mm} 
\end{remark}

\begin{remark}
The Hermite methods studied in this work can be viewed as a technique to improve the accuracy of the primal variable computations 
with mixed element $RT_0$ without resorting to post-processing (see e.g. \cite{BurmanStamm}, \cite{LovadinaStenberg}) or 
hybridization (see e.g. \cite{brezzi} for the diffusion equation and \cite{radu2010} for convection-diffusion problems). 
Incidentally the method proposed in \cite{radu2010} can be applied also to BDM elements to obtain optimal estimates 
\cite{brunner2012} improving  in this way the classical BDM method for convection-diffusion equations \cite{demlow} or for 
stabilization purposes \cite{radu2010, brunner2014} like in several previous work on the subject. Notice that our method allows 
to achieve better accuracy directly from the numerical solution procedure, at negligible additional cost. Thus it seems worthwhile 
searching for Hermite analogs of BDM methods as well in the future. \rule{2mm}{2mm}
\end{remark} 

\noindent \underline{Acknowledgment:}
The first author gratefully acknowledges the support of Statoil through the Akademia agreement, and the second 
author is thankful for the financial support provided by CNPq through grant 307996/2008-5.



\end{document}